\documentclass[11pt]{amsart}
\usepackage{amssymb,amsmath,hyperref}
\usepackage{tikz}

\newtheorem{theorem}{Theorem}[section]
\newtheorem{proposition}[theorem]{Proposition}
\newtheorem{definition}[theorem]{Definition}
\newtheorem{remark}[theorem]{Remark}
\newtheorem{lemma}[theorem]{Lemma}
\newtheorem{corollary}[theorem]{Corollary}
\newtheorem{example}[theorem]{Example}

\newtheorem{procedure}[theorem]{Procedure}

\numberwithin{equation}{section}

\newcommand{\Ppp}{{\mathbb P}}

\newcommand{\binomial}[2]{\binom{#1}{#2}}

\newcommand{\hz}{\widehat{0}}
\newcommand{\ho}{\widehat{1}}
\newcommand{\coveredby}{\prec}

\newcommand{\ptree}{{\mathcal T}}

\DeclareMathOperator{\wt}{wt}
\DeclareMathOperator{\NC}{NC}

\makeatletter
\@namedef{subjclassname@2020}{%
  \textup{2020} Mathematics Subject Classification}
\makeatother

\begin{document}

\vspace*{-20 mm}

\title{On the structure of the $d$-indivisible noncrossing partition posets}
\author{Richard Ehrenborg and G\'abor Hetyei}

\address{Department of Mathematics, University of Kentucky,
Lexington, KY 40506-0027.
\hfill\break \tt http://www.math.uky.edu/\~{}jrge/,
  richard.ehrenborg@uky.edu.}

\address{Department of Mathematics and Statistics, UNC Charlotte,
Charlotte NC 28223-0001.
\hfill\break \tt http://webpages.charlotte.edu/ghetyei/,
ghetyei@charlotte.edu.}



\begin{abstract}
We study the poset of $d$-indivisible noncrossing partitions introduced
by M\"uhle, Nadeau and Williams. These are noncrossing partitions such
that each block and each dual block has cardinality~$1$ modulo $d$. 
Generalizing the work of Speicher, we introduce a generating function
approach to reach new enumerative results and recover some known
formulas on the cardinality, the M\"obius function and the rank numbers.
We compute the antipode of the Hopf algebra of $d$-indivisible
noncrossing partition posets. Generalizing work of Stanley, we show that
the bijection between the maximal chains and $d$-parking functions
introduced by M\"uhle, Nadeau and Williams gives rise to an $EL$-labeling. 
Generalizing a construction of Yan, we also introduce $d$-parking trees
which are in bijective correspondence with the maximal chains.  
\end{abstract}

\maketitle

\section{Introduction}

There has been a lot of work on subposets of the partition lattice.
Most notable is the $d$-divisible partition lattice, which has
been studied by
Stanley~\cite{Stanley_exponential_structures},
Calderbank--Hanlon--Robinson~\cite{Calderbank_Hanlon_Robinson}
and
Wachs~\cite{Wachs}.
See also the papers by~\cite{Ehrenborg_Hedmark,Ehrenborg_Jung}
which extend this lattice to more general partition posets.
In the $d$-divisible partition lattice the size of each block of each partition
is divisible by $d$. Note that that this condition is upward closed
in the partition lattice.
Another condition that is natural to impose on partitions is that
each block size is congruent to $1$ modulo $d$. See for instance
the subposet where all the block sizes are odd
in~\cite{Calderbank_Hanlon_Robinson} 
and also Wachs~\cite[Subsection~4.5.2]{Wachs_Poset_topology}
and the references therein.
Note that the cover relations in such posets may be described by listing 
$d+1$ blocks of a partition which have to be merged into a single block 
to obtain a partition covering the original partition.

The noncrossing partition lattice $\NC_{n}$ is a widely studied subposet of the partition lattice
having many applications;
see the references~\cite{Armstrong,Kreweras,Simion,Simion_Ullman,Speicher,Stanley}.
One of its most striking features is that it is self-dual,
this fact was observed by Kreweras~\cite{Kreweras}
and Simion--Ullman~\cite{Simion_Ullman}.
For an example; see Figure~\ref{figure_example_noncrossing_partition}.
In fact, Stanley pointed out that every interval 
in the noncrossing partition lattice
is self-dual~\cite{Stanley}.
Also note Armstrong studied the subposet of noncrossing partitions where each block
size is divisible by $d$; see~\cite{Armstrong}.

In this paper we will consider the subposet $\NC^{d}_{n}$ of
$d$-indivisible noncrossing partitions of the 
noncrossing partition lattice, first defined by M\"uhle, Nadeau and
Williams~\cite{Muhle_Nadeau_Williams}. These noncrossing partitions are
the ones in which each block size is congruent to $1$ modulo $d$
and each block in the dual partition also
has cardinality congruent to $1$ modulo $d$.
That we require this cardinality condition on each partition and
its dual, makes the subposet~$\NC^{d}_{n}$
naturally self-dual.
On the other hand, many enumerative results for the
noncrossing partition lattice extend to the subposet~$\NC^{d}_{n}$.
One less studied invariant is the antipode of a poset
which is the natural Hopf algebra extension of the M\"obius function.
For the poset $\NC^{d}_{n}$ we obtain its antipode
in terms of noncrossing hypertrees,
generalizing the results in~\cite{Ehrenborg_Happ}.
We also enumerate the maximal chains in the poset $\NC^{d}_{n}$
and show that via an edge labeling they are in
bijection with $d$-parking functions, as defined
by Stanley~\cite{Stanley_k-parking} and Yan~\cite{Yan}.
Furthermore, these structures are also in bijective
correspondence with a class of trees that we call $d$-parking trees.

The paper is organized as follows.
In Section~\ref{section_Preliminaries}
we present basic results about the noncrossing partition lattice~$\NC_{n}$
including the Kreweras dual and review basic generating function results.
Our main object of study, the poset $\NC^{d}_{n}$ is
introduced in Section~\ref{section_NC^{d}_{n}}. We review some of
the basic results presented in~\cite{Muhle_Nadeau_Williams}. We prove
that intervals in this poset are products of 
smaller posets of the same form.
In Section~\ref{section_Enumerative}
we obtain enumerative results on~$\NC^{d}_{n}$.
Our approach is using generating functions; 
see Subsection~\ref{subsection_Generating}.
This method generalizes results of Speicher~\cite{Speicher}; 
see Proposition~\ref{proposition_Speicher_d}.
In Section~\ref{section_A_tree_representation_of_noncrossing_partitions}
we revisit a labeled tree representation of noncrossing partitions which
also appears in the work of M\"uhle, Nadeau and
Williams~\cite{Muhle_Nadeau_Williams}, but we use a different recursive
algorithm relying on selecting a root vertex instead of using their root
edge selection method. Our labeled plane tree model
may be used to visualize the proofs of our enumerative results in in
Section~\ref{section_Enumerative}, to this extent our rooted trees play
the same role as the differently rooted trees
in~\cite{Muhle_Nadeau_Williams}. Our root vertex selection method is 
somewhat similar to the ``dual tree representation''
discussed by Kortchemski and Marzouk~\cite[Section 2]{Kortchemski_Marzouk}
relying on a bijection introduced by Janson and
Stef\'ansson~\cite{Janson_Stefansson}, who use their recursive
construction to build an algorithm that generates a random noncrossing
partition in a uniform fashion.  
In Section~\ref{section_antipode}, 
we obtain an expression for the antipode of $\NC^{d}_{n}$
in terms of noncrossing hypertrees.
In Section~\ref{section_EL-labeling_of_the_poset}
we show that the poset~$\NC^{d}_{n}$
has an $EL$-labeling.
This labeling has a connection with a bijection
between maximal chains and
$d$-parking functions,
going back to the papers~\cite{Muhle_Nadeau_Williams}
and Irving and Rattan~\cite{Irving_Rattan}.
We introduce the notion of a $d$-parking tree in
Section~\ref{section_tree_representation}. These labeled trees are in bijection
with the $d$-parking functions and visually encode the block merging
processes represented by the maximal chains in $\NC^{d}_{n}$.
Finally, further research and open problems are
outlined in Section~\ref{section_Concluding_remarks}.

\section{Preliminaries}
\label{section_Preliminaries}

\subsection{The noncrossing partition lattice}
\label{subsection_The_noncrossing_partition_lattice}

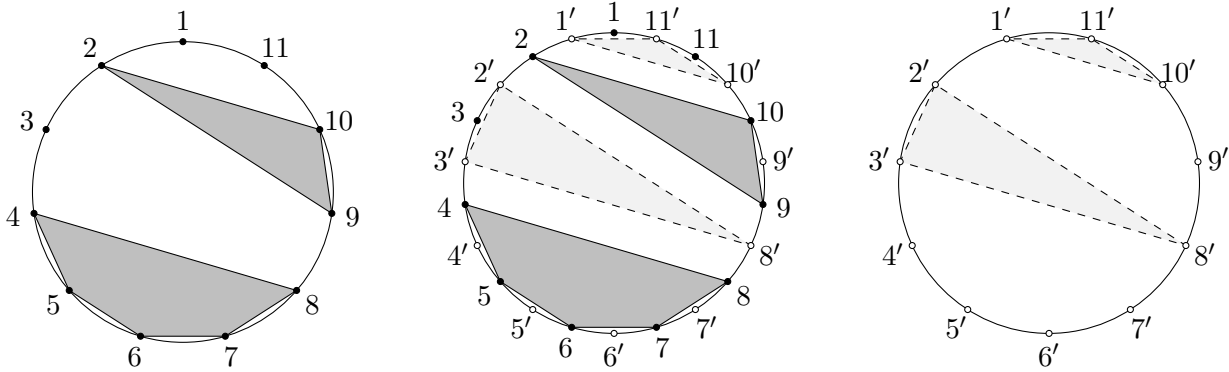
\begin{figure}
\newcommand{\app}{360/11}
\newcommand{\bpp}{90-360/11}
\newcommand{\gpp}[1]{({5*cos(\app*#1+\bpp)},{5*sin(\app*#1+\bpp)})}
\newcommand{\hpp}[1]{({5.7*cos(\app*#1+\bpp)},{5.7*sin(\app*#1+\bpp)})}
\newcommand{\bqq}{90-360/11+360/22}
\newcommand{\gqq}[1]{({5*cos(\app*#1+\bqq)},{5*sin(\app*#1+\bqq)})}
\newcommand{\hqq}[1]{({5.7*cos(\app*#1+\bqq)},{5.7*sin(\app*#1+\bqq)})}
\newcommand{\rrr}{0.24}
\begin{center}
\begin{tikzpicture}[scale = 0.40,inner sep=0.1mm]
\draw[-] (0,0) circle (5);

\draw[fill=gray!50] \gpp{2} -- \gpp{9} -- \gpp{10} -- cycle;
\draw[fill=gray!50] \gpp{4} -- \gpp{5} -- \gpp{6} -- \gpp{7} -- \gpp{8} -- cycle;

\foreach \k in {1, 2, ..., 11}
{
\node at \hpp{\k} {$\k$};
\draw[-,fill] \gpp{\k} circle (0.1);
};
\end{tikzpicture}
\hspace*{7 mm}
\begin{tikzpicture}[scale = 0.40,inner sep=0.1mm]
\draw[-] (0,0) circle (5);

\draw[fill=gray!50] \gpp{2} -- \gpp{9} -- \gpp{10} -- cycle;
\draw[fill=gray!50] \gpp{4} -- \gpp{5} -- \gpp{6} -- \gpp{7} -- \gpp{8} -- cycle;
\draw[dashed,fill=gray!10] \gqq{1} -- \gqq{10} -- \gqq{11} -- cycle;
\draw[dashed,fill=gray!10] \gqq{2} -- \gqq{3} -- \gqq{8} -- cycle;

\foreach \k in {1, 2, ..., 11}
{
\node at \hpp{\k} {$\k$};
\node at \hqq{\k} {$\k^{\prime}$};
\draw[-,fill] \gpp{\k} circle (0.1);
\draw[-,fill=white] \gqq{\k} circle (0.1);
};
\end{tikzpicture}
\hspace*{7 mm}
\begin{tikzpicture}[scale = 0.40,inner sep=0.1mm]
\draw[-] (0,0) circle (5);

\draw[dashed,fill=gray!10] \gqq{1} -- \gqq{10} -- \gqq{11} -- cycle;
\draw[dashed,fill=gray!10] \gqq{2} -- \gqq{3} -- \gqq{8} -- cycle;

\foreach \k in {1, 2, ..., 11}
{
\node at \hqq{\k} {$\k^{\prime}$};
\draw[-,fill=white] \gqq{\k} circle (0.1);
};
\end{tikzpicture}
\end{center}
\caption{The noncrossing partition $\pi$ given by
$1|2,9,10|3|4,5,6,7,8|11$,
its circle representation
and the Kreweras dual partition
$1^{\prime},10^{\prime},11^{\prime}|2^{\prime},3^{\prime},8^{\prime}|4^{\prime}|5^{\prime}|6^{\prime}|7^{\prime}|9^{\prime}$.
}
\label{figure_example_noncrossing_partition}
\end{figure}

The partition lattice consists of all partitions of
the set $[n] = \{1,2, \ldots, n\}$.
Two different blocks $B$ and~$C$ are noncrossing if there are no four
elements $i < j < k < \ell$ such that $i,k \in B$ and $j,\ell \in C$.
This condition is viewed geometrically as follows. Place the elements of the set $[n]$
in a circle in positive orientation. The two blocks $B$ and~$C$ are noncrossing
if the convex hull of the two blocks are disjoint.
A partition is noncrossing if all pairs of blocks are noncrossing.
For an example,
the partition
$\pi = \{\{1\}, \{2,9,10\}, \{3\}, \{4,5,6,7,8\}, \{11\}\}$
is displayed in the first drawing in
Figure~\ref{figure_example_noncrossing_partition}.
For short hand we write this partition as
$1|2,9,10|3|4,5,6,7,8|11$.
We call these drawings of the noncrossing partition and its dual
the {\em circle representation} of the partition.

The {\em Kreweras dual} of a noncrossing partition is defined as follows.
First observe that the complement of the union of the convex hull of the blocks
form regions. For instance,
in the partition~$\pi$ in
Figure~\ref{figure_example_noncrossing_partition}
we have a region $R$ with the vertices $2$, $3$, $4$, $8$ and $9$.
However, to obtain the dual partition place the vertex $i^{\prime}$
between vertices $i$ and $i+1$.
That is, we have the cyclic order
$1 < 1^{\prime} <
2 < 2^{\prime} <
3 < 3^{\prime} <
\cdots <
n-1 < (n-1)^{\prime} <
n < n^{\prime} < 1$.
The region $R$ now correspond to the block
$\{2^{\prime},3^{\prime},8^{\prime}\}$.
The blocks corresponding to the regions form the dual partition
which is noncrossing. In our example, we obtain
the dual partition $\pi^{\prime}$ given by
$1^{\prime},10^{\prime},11^{\prime}|2^{\prime},3^{\prime},8^{\prime}|4^{\prime}|5^{\prime}|6^{\prime}|7^{\prime}|9^{\prime}$.

We say that a block $B$ of the noncrossing partition $\pi$ and
a block $C^{\prime}$ of the dual partition are {\em adjacent}
if the block $B$ contains the two elements $i$ and $j$
(they could be the same element)
and the block $C^{\prime}$
contains the elements
$i^{\prime}$
and
$(j-1)^{\prime}$.
Note that the $4$-gon
with vertices
$i$, $i^{\prime}$, $(j-1)^{\prime}$ and $j$
forms a channel between the block~$B$ and the dual block~$C^{\prime}$.

Following~\cite[Definitions~3.1 and~3.2]{Ehrenborg_Happ}
we make the following definition.
For a block $B$ in a noncrossing partition $\pi$ and
an adjacent block $C^{\prime}$ in the dual partition $\pi^{\prime}$
let $i$ and $j$ be the two vertices of the block~$B$
that are adjacent to the block~$C^{\prime}$
and pick
$\gamma(B,C^{\prime})$ to be the vertex
of these two that is the most
clockwise orientation from the block $B$'s perspective.

A small note is in order at this point.
The notion related to $\gamma(B,C^{\prime})$ introduced
in~\cite{Ehrenborg_Happ} is denoted $v_{R}(B)$. The
definition of $\gamma(B,C^{\prime})$  differs in two key points. 
First, here we work with the blocks of the dual partition,
whereas~\cite{Ehrenborg_Happ} uses regions in the partition. Second,
in~\cite{Ehrenborg_Happ} the authors pick the vertex
in the most positive orientation. Here we pick the vertex in the most
negative orientation. The reason for this difference is that
we would like to be consistent with the Kreweras notation of the dual;
see Section~\ref{section_A_tree_representation_of_noncrossing_partitions}.

Let $|\pi|$ denote the number of blocks of the partition $\pi$.
Let $\NC_{n}$ be the collection all noncrossing partitions on the set $[n]$
ordered by refinement. This poset is graded and is in fact a self-dual lattice.
That is, the lattice $\NC_{n}$ is self-dual.
Also note that rank of a noncrossing partition $\pi$ is given by
$\rho(\pi) = n - |\pi|$
and that the identity $|\pi| + |\pi'| = n+1$ holds.

Consider a pair of noncrossing partitions $\pi$ and $\sigma$ in the
noncrossing partition lattice~$\NC_{n}$ such that $\pi \leq \sigma$.
The interval $[\pi,\sigma]$ is given by $\{\tau \in \NC_{n} : \pi \leq \tau \leq \sigma\}$.
It is well-known that intervals in the noncrossing partition lattice are isomorphic to
Cartesian products of noncrossing partition lattices~\cite[Proposition~1]{Speicher}.

\subsection{Generating functions}

We need the following classical results about generating functions.
Assume that the generating function $G(z)$ satisfies the algebraic equation
\begin{align}
G(z) & = z \cdot G(z)^{m}   + 1
\label{equation_algebraic_equation}
\end{align}
then the coefficient of $z^{n}$ of $G(z)^{\ell}$
is given by
\begin{align}
[z^{n}] G(z)^{\ell}
& =
\frac{\ell}{mn+\ell} \cdot \binom{mn+\ell}{n}
=
\frac{\ell}{mn+\ell-n} \cdot \binom{mn+\ell-1}{n} ,
\label{equation_coefficient_of_the_generating_function}
\end{align}
for $\ell \geq 1$;
see~\cite[Equations~(7.68) and~(7.69)]{Graham_Knuth_Patashnik}.
See also~\cite[Example~6.2.6]{Stanley_EC2}.

The Lagrange inversion formula allows one to obtain
the coefficients of a power series $f(z)$ that satisfies
the functional equation $f(z) = z \cdot g(f(z))$.
Good stated a multivariate extension of the Lagrange inversion formula~\cite{Good}.
See also~\cite[Theorem~4]{Ehrenborg_Mendez}.
We will only need this extension in two variables and
hence we state it only in this case.
\begin{theorem}[Good's Inversion Formula for two variables.]
Let $g_{1}(z,w)$ and $g_{2}(z,w)$ be formal power series
with constant coefficients.
Assume that the formal power series $f_{1}(z,w)$ and $f_{2}(z,w)$ satisfy the equations
\begin{align*}
f_{1}(z,w)
& =
z \cdot g_{1}(f_{1}(z,w),f_{2}(z,w)) ,
&
f_{2}(z,w)
& =
w \cdot g_{2}(f_{1}(z,w),f_{2}(z,w)) .
\end{align*}
Then the coefficient of $z^{m} w^{n}$ of $f_{1}(z,w)^{k} \cdot f_{2}(z,w)^{\ell}$
is given by
\begin{align*}
[z^{m} w^{n}] f_{1}(z,w)^{k} \cdot f_{2}(z,w)^{\ell}
& =
[z^{m-k} w^{n-\ell}]
\det
\begin{pmatrix}
g_{1}^{m} - z \cdot \frac{\partial g_{1}}{\partial z} \cdot g_{1}^{m-1}
&
- w \cdot \frac{\partial g_{1}}{\partial w} \cdot g_{1}^{m-1}
\\
- z \cdot \frac{\partial g_{2}}{\partial z} \cdot g_{2}^{n-1}
&
g_{2}^{n} - w \cdot \frac{\partial g_{2}}{\partial w} \cdot g_{2}^{n-1}
\end{pmatrix} .
\end{align*}
\label{theorem_Good}
\end{theorem}

\section{The poset of $d$-indivisible noncrossing partitions}
\label{section_NC^{d}_{n}}

The following concept was introduced by M\"uhle, Nadeau and
Williams~\cite{Muhle_Nadeau_Williams}.  

\begin{definition}
  Let $d$ be a positive integer. A noncrossing partition $\pi\in\NC_{n}$
  is {\em $d$-indivisible} if all the blocks of the partition $\pi$ and
  the dual partition $\pi^{\prime}$ have cardinality congruent to $1$
  modulo $d$. We denote the subposet of $d$-indivisible noncrossing
  partitions of $\NC_{n}$ by $\NC_{n}^{d}$.
\end{definition}  
As an example, the partition in
Figure~\ref{figure_example_noncrossing_partition}
belongs the poset $\NC^{2}_{11}$
since all the blocks in the partition and the dual partition
have odd cardinality.

In~\cite{Muhle_Nadeau_Williams} $d$-indivisible noncrossing partitions are only defined
in the case when $n \equiv 1 \bmod d$ holds. The next lemma shows
that a generalization to other values of $n$ is not possible. 
\begin{lemma}
If the set $\NC^{d}_{n}$ is nonempty then $n \equiv 1 \bmod d$.
\label{lemma_fundamental}
\end{lemma}
\begin{proof}
Assuming that $\NC^{d}_{n}$ is nonempty,
let $\pi$ be a noncrossing partition in $\NC^{d}_{n}$.
Then we have the following string of congruences
$n = \sum_{B \in \pi} |B| \equiv  \sum_{B \in \pi} 1 = |\pi| \bmod d$.
Similarly, we obtain
$n \equiv |\pi'| \bmod d$.
Adding these two congruences yields
$2n \equiv |\pi| + |\pi'| = n+1 \bmod d$, that is, 
$n \equiv 1 \bmod d$.
\end{proof}

Note that when $n \equiv 1 \bmod d$ then $\NC^{d}_{n}$
is nonempty since it contains the noncrossing partition~$\{[n]\}$.
Furthermore, this partition~$\{[n]\}$ is the maximal element in $\NC^{d}_{n}$
and the noncrossing partition $\{\{1\}, \{2\}, \ldots,\{n\}\}$ is the
minimal element.

By the same the same reasoning as in Lemma~\ref{lemma_fundamental}
we have the next lemma.
\begin{lemma}
Let $n$ be a positive integer such that $n \equiv 1 \bmod d$ and
let $\pi$ be a noncrossing partition in~$\NC_{n}$.
Assume that all blocks $B$ in $\pi$ satisfy $|B| \equiv 1 \bmod d$.
Assume furthermore
that all blocks~$C'$ but possibly one in the dual partition $\pi'$ satisfy $|C'| \equiv 1 \bmod d$.
Then the final block also has size
congruent to $1 \bmod d$ and
the noncrossing partition~$\pi$ belongs to~$\NC^{d}_{n}$.
\label{lemma_all_but_one}
\end{lemma}
\begin{proof}
Suppose that the dual partition $\pi'$ contains exactly
one block $D'$ such that
$|D'| \not\equiv 1 \bmod d$.
We obtain a contradiction 
by the following string of congruences:
\begin{align*}
|D'|
& =
n - \sum_{\substack{C' \in \pi' \\ C' \neq D'}} |C'|
\equiv
n - (|\pi'| - 1)
=
n + 1 - |\pi'|
\equiv
|\pi|
\equiv
\sum_{B \in \pi} |B|
\equiv 
n
\equiv
1 \bmod d.
\qedhere
\end{align*}
\end{proof}

\begin{lemma}
Assume that the noncrossing partition $\pi$ belongs to the noncrossing partition poset $\NC^{d}_{n}$.
Then the two following equalities hold:
\begin{align*}
\frac{|\pi'| - 1}{d} 
& =
\sum_{B \in \pi} \frac{|B| - 1}{d} ,
&
\frac{|\pi| - 1}{d} 
& =
\sum_{C' \in \pi'} \frac{|C'| - 1}{d} ,
\end{align*}
where all the summands are integers.
\label{lemma_blocks_versus_dual_partition}
\end{lemma}
\begin{proof}
By duality it is enough to prove the first identity.
We have $|\pi'| - 1 = n - |\pi| = \sum_{B \in \pi} (|B| - 1)$,
and the result follows by dividing by $d$.
By the definition of $\NC^{d}_{n}$, all summands in the statement are integers.
\end{proof}

From~\cite{Muhle_Nadeau_Williams} we have:
\begin{theorem}[M\"uhle--Nadeau--Williams]
The noncrossing partition poset $\NC^{d}_{n}$ is graded of rank $(n-1)/d$
and the rank function is given by $\rho(\pi) = (n-|\pi|)/d$.
Hence for a cover relation $\pi \coveredby \sigma$
the partition $\sigma$ is obtained by joining $d+1$ blocks of $\pi$.
Furthermore, the Kreweras dual
is an order reversing bijection on $\NC^{d}_{n}$.
\label{theorem_NC_structure}
\end{theorem}

\begin{corollary}
If $\sigma\in\NC^{d}_{n}$ is a coatom then it consists of $d+1$ blocks
and these blocks are all intervals in the cyclic order on $1,2, \ldots, n$.
\label{corollary_coatom}
\end{corollary}
This corollary follows from~\cite[Corollary~2.3]{Muhle_Nadeau_Williams}.
A direct argument is as follows.
\begin{proof}[Proof of Corollary~\ref{corollary_coatom}.]
Since an atom in the poset consists of one non-singleton block of cardinality $d+1$
and the remaining blocks are singletons, the result follows by taking the dual partition.
\end{proof}

\begin{lemma}
Let  $\sigma$ be a coatom in the noncrossing partition poset $\NC^{d}_{n}$,
and assume that $\sigma = \{B_{1},B_{2}, \ldots, B_{d+1}\}$.
Then the interval $[\hz,\sigma]$ in $\NC^{d}_{n}$ factors as
\begin{align*}
[\hz,\sigma]
& \cong
\prod_{i=1}^{d+1} \NC^{d}_{|B_{i}|} .
\end{align*}
\label{lemma_one_small_step}
\end{lemma}
\begin{proof}
Let $D^{\prime}$ be the unique singleton block of the dual partition $\sigma^{\prime}$.
Let $\tau$ be a noncrossing partition in the interval~$[\hz,\sigma]$,
that is, $\tau \leq \sigma$.
Let $\tau_{i}$ be the restriction of the partition $\tau$ to the block $B_{i}$.
Note that $|B_{i}| \equiv 1 \bmod d$.
Furthermore, each block $B$ of $\tau_{i}$
is a block of $\tau$ and hence satisfies $|B| \equiv 1 \bmod d$.
Let $\nu_{i}$ be the partition
$\tau_{i} \cup \{B_{1}, \ldots, \widehat{B_{i}}, \ldots, B_{d+1}\}$.
Directly the inequalities $\tau \leq \nu_{i} \leq \sigma$ hold.
Hence the dual partition $\nu^{\prime}_{i}$ has a unique block $C^{\prime}_{i}$ that contains
the coblock $D^{\prime}$.
Restricting the partition $\nu^{\prime}_{i}$ to $\tau^{\prime}_{i}$
we observe that all but one of the blocks of the dual partition $\tau^{\prime}_{i}$
is a block of the dual partition $\nu^{\prime}_{i}$ and hence of $\tau^{\prime}$.
The only block left out is the restriction of~$C^{\prime}_{i}$.
Hence we know that all but one of the dual blocks~$C^{\prime}$ of $\tau_{i}$ satisfies  
$|C^{\prime}| \equiv 1 \bmod d$.
By Lemma~\ref{lemma_all_but_one}
we conclude that $\tau_{i}$ belongs to $\NC^{d}_{|B_{i}|}$.

Furthermore, the map from the interval~$[\hz,\sigma]$
to the Cartesian product $\prod_{i=1}^{d+1} \NC^{d}_{|B_{i}|}$,
sending $\tau$ to the $(d+1)$-tuple
$(\tau_{1},\tau_{2}, \ldots, \tau_{d+1})$ is bijective
and order preserving.
\end{proof}

A key result about the noncrossing partition poset is
the following theorem.
\begin{theorem}
Let $\pi$ and $\sigma$ be two noncrossing partitions in $\NC^{d}_{n}$
such that $\pi \leq \sigma$ holds.
Then the interval $[\pi,\sigma]$ in $\NC^{d}_{n}$ is isomorphic to
the Cartesian product of smaller
posets of the form $\NC^{d}_{k}$.
\label{theorem_Cartesian_product}
\end{theorem}
\begin{proof}
The result follows directly from
Lemma~\ref{lemma_one_small_step}
and the dual version of Lemma~\ref{lemma_one_small_step}.
\end{proof}

\begin{remark}
{\rm For two sets $B \subseteq \{1,2,  \ldots, n\}$
and $C' \subseteq \{1',2',  \ldots, n'\}$
define the {\em intertwining number} $i(B,C')$
to be the number of times we move from the set $B$ to set $C'$
when reading the entries of $B \cup C'$ in the
cyclic order
$1 < 1^{\prime} <
2 < 2^{\prime} <
3 < 3^{\prime} <
\cdots <
n-1 < (n-1)^{\prime} <
n < n^{\prime} < 1$.
Note that $i(B,C') = 0$ if and only if $B$ or $C'$ is the empty set.
Let $\pi$ and $\sigma$ be two noncrossing partitions in $\NC^{d}_{n}$
such that $\pi \leq \sigma$.
Then the interval $[\pi,\sigma]$ in $\NC^{d}_{n}$ is isomorphic to
the Cartesian product
\begin{align*}
[\pi,\sigma]
& \cong
\prod_{B \in \sigma}
\prod_{C' \in \pi'}
\NC^{d}_{i(B,C')} .
\end{align*}
However, we do not need this statement and hence omit its proof.}
\end{remark}

The following characterization of $d$-indivisible noncrossing
partition is part of~\cite[Theorem~1.1]{Muhle_Nadeau_Williams}.

\begin{theorem}[M\"uhle--Nadeau--Williams]
\label{theorem_characterization}
The noncrossing partition $\pi$ belongs to $\NC^{d}_{n}$
if and only if
the following two conditions are satisfied
(i)
$n \equiv 1 \bmod d$;
(ii)
for each pair of cyclically consecutive elements $i$ and $j$
of a nonsingleton block $B$ in a noncrossing partition $\pi$,
the number of elements in $[n]$ strictly between $i$ and $j$ in the cyclic order
is divisible by~$d$.
\end{theorem}

In~\cite[Theorem~1.1]{Muhle_Nadeau_Williams} a third equivalent characterization was
added in terms of permutation factorizations. We note that 
Theorem~\ref{theorem_characterization}
can be proved directly using the results in the section.

\section{Enumerative results for the poset $NC^{d}_{n}$}
\label{section_Enumerative}

\subsection{Generating function approach}
\label{subsection_Generating}

We now introduce generating functions in order to obtain
enumerative results for the noncrossing partition poset $\NC^{d}_{n}$.
Let $a(0) = 1, a(1), a(2), \ldots$
and $a^{*}(0) = 1, a^{*}(1), a^{*}(2), \ldots$
be two given sequences.
For a noncrossing partition $\pi$ in the poset $\NC^{d}_{dk+1}$
define its weight by the product:
\begin{align*}
\wt(\pi)
& =
\prod_{B \in \pi} a\left(\frac{|B|-1}{d}\right)
\cdot
\prod_{C^{\prime} \in \pi^{\prime}} a^{*}\left(\frac{|C^{\prime}|-1}{d}\right)
\cdot
s^{({|\pi|-1})/{d}}
\cdot
t^{({|\pi^{\prime}|-1})/{d}} .
\end{align*}
Note that the power of the variable $s$ is the corank of the partition $\pi$
and that the power of $t$ is the rank of $\pi$
in the poset $\NC^{d}_{dk+1}$.
We are interested in  exploring the following sum and its ordinary generating function:
\begin{align}
b(k)
& =
\sum_{\pi \in \NC^{d}_{dk+1}}
\wt(\pi) ,
&
B(x)
& =
\sum_{k \geq 0} b(k) \cdot x^{k} .
\label{equation_b_k}
\end{align}
Hence every monomial $s^{i} t^{j} \cdot x^{k}$ that occurs
in the generating function $B(x)$ satisfies $i+j=k$.
Let $A(x)$ and $A^{*}(x)$ be the two ordinary generating functions
\begin{align*}
A(x) & = \sum_{k \geq 0} a(k) \cdot x^{k} ,
&
A^{*}(x) & = \sum_{k \geq 0} a^{*}(k) \cdot x^{k} .
\end{align*}
In order to obtain information about $B(x)$, we consider the two following partial sums
and their generating functions
\begin{align*}
c(k)
& =
\sum_{\substack{\pi \in \NC^{d}_{dk+1} \\ \{1\} \in \pi}}
\wt(\pi) ,
&
C(x)
& =
\sum_{k \geq 0} c(k) \cdot x^{k} , \\
c^{*}(k)
& =
\sum_{\substack{\pi \in \NC^{d}_{dk+1} \\ \{1^{\prime}\} \in \pi^{\prime}}}
\wt(\pi) ,
&
C^{*}(x)
& =
\sum_{k \geq 0} c^{*}(k) \cdot x^{k} .
\end{align*}
That is, $c(k)$ and its generating function $C(x)$, enumerate
noncrossing partitions that have a given element as a singleton block.
Similarly, $c^{*}(k)$ and $C^{*}(x)$ enumerate partitions
with a given element as a singleton block in the dual.

\begin{theorem}
The following three generating function identities hold:
\begin{align}
\label{equation_C_and_C*}
C(x)
& =
A^{*}\left(x \cdot s \cdot C^{*}(x)^{d}\right) ,
&
C^{*}(x)
& =
A\left(x \cdot t \cdot C(x)^{d}\right) , \\
B(x) & = C(x) \cdot C^{*}(x) .
\end{align}
Furthermore, the two identities in~\eqref{equation_C_and_C*} uniquely determine
the generating functions $C(x)$ and $C^{*}(x)$.
\label{theorem_generating_functions_C_and_C*}
\end{theorem}
\begin{proof}
Observe that the first two identities are dual to each other. Hence it
is enough 
that we prove the second identity.
Let $\pi$ be a partition in $\NC^{d}_{dk+1}$ where $dk+1 = n$
and where the singleton block~$\{n^{\prime}\}$
belongs to the dual partition $\pi^{\prime}$.
That is, the partition contains a block $B$ such that
\begin{align*}
B
& =
\{1 = i_{1} < i_{2} < \cdots < i_{dr+1} = n\} .
\end{align*}
By Theorem~\ref{theorem_characterization}
each difference $i_{j+1} - i_{j}$ is congruent to $1$ modulo $d$.
Let $\tau_{j}$ now be the restriction of the partition $\pi$
to the union $[i_{j},i_{j+1}] \cup B$ and identify all the elements of the block $B$.
This construction is equivalent to letting $\tau_{j}$ be the restriction to the
open interval $(i_{j},i_{j+1})$ and adding a singleton block.
Keeping track of the blocks in~$\pi$ we have that
\begin{align*}
|\pi| - 1
& =
(|\tau_{1}| - 1)
+ \cdots +
(|\tau_{dr}| - 1) .
\end{align*}
Let $\tau_{i}^{\prime}$ be the dual partition of $\tau_{i}$ inside the set $[i_{j},i_{j+1}] \cup B$.
Counting the blocks in the dual partition we obtain
\begin{align*}
|\pi^{\prime}| - 1
& =
|\tau^{\prime}_{1}|
+ \cdots +
|\tau^{\prime}_{dr}|
=
dr
+
(|\tau^{\prime}_{1}| - 1)
+ \cdots +
(|\tau^{\prime}_{dr}| - 1) .
\end{align*}
Hence the weight $\wt(\pi)$ factors as
\begin{align*}
\wt(\pi)
& =
a(r)
\cdot
t^{r}
\cdot
\wt(\tau_{1})
\cdots
\wt(\tau_{dr}) .
\end{align*}
Now we express $c^{*}(k)$ as the sum
\begin{align}
\label{equation_recursion_for_c*}
c^{*}(k)
& =
\sum_{r \geq 1}
a(r)
\cdot
t^{r}
\cdot
\sum_{
\substack{1 = i_{1} < i_{2} < \cdots < i_{dr+1} = n \\
i_{j+1} - i_{j} \equiv 1 \bmod d}
}
\sum_{\tau_{1}}
\cdots
\sum_{\tau_{dr}}
\prod_{j=1}^{dr}
\wt(\tau_{j}) \\
\nonumber
& =
\sum_{r \geq 1}
a(r)
\cdot
t^{r}
\cdot
\sum_{
\substack{1 = i_{1} < i_{2} < \cdots < i_{dr+1} = n \\
i_{j+1} - i_{j} \equiv 1 \bmod d}
}
\prod_{j=1}^{dr}
c\left( \frac{i_{j+1} - i_{j} - 1}{d} \right) ,
\end{align}
where each $\tau_{j}$ ranges over partitions in 
$\NC^{d}_{i_{j+1}-i_{j}}$
with a given singleton block.
Observe that $k$ is given by the sum
\begin{align*}
k
& =
r + \frac{i_{2} - i_{1} - 1}{d}
+ \cdots + \frac{i_{dr+1} - i_{dr} - 1}{d} .
\end{align*}
Hence multiply with $x^{k}$ and sum over all $k \geq 1$.
\begin{align*}
\sum_{k \geq 1} c^{*}(k) \cdot x^{k}
& =
\sum_{r \geq 1}
a(r) \cdot t^{r} \cdot x^{r} \cdot
\sum_{
\substack{1 = i_{1} < i_{2} < \cdots < i_{dr+1} = n \\
i_{j+1} - i_{j} \equiv 1 \bmod d}
}
\prod_{j=1}^{dr}
c\left( \frac{i_{j+1} - i_{j} - 1}{d} \right) 
\cdot
x^{({i_{j+1} - i_{j} - 1})/{d}} \\
& =
\sum_{r \geq 1}
a(r) \cdot t^{r} \cdot x^{r} \cdot C(x)^{dr} .
\end{align*}
Lastly, adding the constant $c^{*}(0) = 1$ yields the second equation.

Equation~\eqref{equation_recursion_for_c*} and its dual
yield a pair of intertwined recursions for the sequences
$(c(k))_{k\geq 0}$ and $(c^{*}(k))_{k\geq 0}$: 
equation~\eqref{equation_recursion_for_c*}
expresses $c(k)$ 
in terms of $c^{*}(0), \ldots, c^{*}(k-1)$, whereas the 
dual equation
expresses $c^{*}(k)$
in terms of $c(0), \ldots, c(k-1)$.
Together with the initial condition $c(0) = c^{*}(0) = 1$,
the two functions $C(x)$ and $C^{*}(x)$ are uniquely determined.

Now we prove the third identity.
Given a partition $\pi$ in $\NC^{d}_{dk+1}$,
let $j$ be the largest element in the block $B$ containing $1$.
Theorem~\ref{theorem_characterization} implies
that $j \equiv 1 \bmod d$, hence let $j = dr+1$.
Similarly, then $j^{\prime}$ is the smallest element in
block of the dual partition containing $n^{\prime}$.
Let $\sigma$ be the restriction $\pi$ to the interval $[1,j]$,
which has cardinality $dr+1$.
Let $\tau$ be the restriction of $\pi$ to the union $[j+1,n] \cup B$
and identify all elements of the block~$B$.
Note that $\tau$ is a partition on a set of size $d(k-r) + 1$.
Note that $\sigma$ has a given singleton block in the dual
and $\tau$ has also a given singleton block.
Enumerating the blocks in $\pi$ and the blocks in the dual partition
we have
\begin{align*}
|\pi| -1 & = |\sigma| -1 + |\tau| - 1 , &
|\pi^{\prime}| -1 & = |\sigma^{\prime}| -1 + |\tau^{\prime}| - 1 .
\end{align*}
Hence the weight of $\pi$ factors as
\begin{align*}
\wt(\pi) 
& =
\wt(\tau) \cdot \wt(\sigma) .
\end{align*}
Summing over all partitions $\pi$ in $\NC^{d}_{dk+1}$, multiplying with
$x^{k} = x^{r} \cdot x^{k-r}$ and summing over all
$k \geq 0$ yield the desired identity.
\end{proof}

\begin{proposition}
Set $s = t = 1$ and $a^{*}(k) = (-1)^{k} \cdot \binom{-1/d}{k}$ for all $k$.
Define the two functions
$\overline{A}(x) = 1 + x \cdot A(x)^{d}$
and
$\overline{B}(x) = 1 + x \cdot B(x)^{d}$.
Then the following identity holds
\begin{align*}
\overline{B}(x) & = \overline{A}(x \cdot \overline{B}(x)) .
\end{align*}
\label{proposition_Speicher_d}
\end{proposition}
In the case when $d=1$, this is Speicher's identity;
see~\cite[Theorem on page~616]{Speicher}.
See also earlier work of Beissinger~\cite[Theorem~6.7]{Beissinger}.
\begin{proof}[Proof of Proposition~\ref{proposition_Speicher_d}]
Note that $A^{*}(x) = 1/\sqrt[d]{1-x}$,
that is,
$A^{*}(x)^{d} = 1/(1-x)$.
Hence we have
$C(x)^{d} = 1/(1 - x \cdot C^{*}(x)^{d})$.
Cross multiply to obtain
$1 = C(x)^{d} - x \cdot C(x)^{d} \cdot C^{*}(x)^{d} = C(x)^{d} - x \cdot B(x)^{d}$.
We obtain
$\overline{B}(x) = 1 + x \cdot B(x)^{d} = C(x)^{d}$.
Now we have
\begin{align*}
\overline{B}(x)
& =
1 + x \cdot B(x)^{d}
=
1 + x \cdot C(x)^{d} \cdot C^{*}(x)^{d}
=
1 + x \cdot C(x)^{d} \cdot A(x \cdot C(x)^{d})^{d} \\
& =
1 + x \cdot \overline{B}(x) \cdot A(x \cdot \overline{B}(x))^{d}
=
\overline{A}(x \cdot \overline{B}(x)) .
\qedhere
\end{align*}
\end{proof}

\begin{example}
{\rm
The following results of M\"uhle, Nadeau and
Williams~\cite{Muhle_Nadeau_Williams}     
follow from Theorem~\ref{theorem_generating_functions_C_and_C*}.
\begin{itemize}
\item[(a)]
\cite[Theorem~1.2]{Muhle_Nadeau_Williams}
The cardinality of $\NC^{d}_{dk+1}$ is given by
\begin{align*}
\left| \NC^{d}_{dk+1} \right|
& =
\frac{2}{dk+2} \cdot \binom{dk+k+1}{k} .
\end{align*}
This follows by setting $a(k) = a^{*}(k) = s = t = 1$
and using the relations
$A(x) = A^{*}(x) = 1/(1-x)$
and $C(x) = C^{*}(x)= 1 + x \cdot C(x)^{d+1}$.
\item[(b)]
\cite[Corollary~4.6]{Muhle_Nadeau_Williams}
The number of partitions of rank $j$ 
in the noncrossing partition poset $\NC^{d}_{dk+1}$,
for $k = i+j$, is given by
\begin{align*}
\frac{dk+1}{(i+dj+1) \cdot (di+j+1)}
\cdot
\binomial{i+dj+1}{i} \cdot \binomial{di+j+1}{j} .
\end{align*}
Set $a(k) = a^{*}(k) = 1$ such that
$A(x) = A^{*}(x) = 1/(1-x)$.
Let $D(s,t,x) = C(s,t,x) - 1$ and $D^{*}(s,t,x) = C^{*}(s,t,x) - 1$.
Rewrite the two equations in~\eqref{equation_C_and_C*} as
\begin{align*}
D
& =
sx \cdot (1+D) \cdot (1+D^{*})^{d} , &
D^{*}
& =
tx \cdot (1+D)^{d} \cdot (1+D^{*}) .
\end{align*}
Note that every monomial that occurs has the form $s^{i} t^{j} \cdot x^{i+j}$.
At this point we set $x=1$ since all the information is carried by the coefficients
of the monomials in the variables $s$ and~$t$.
We can now determine the coefficients of $(1+D) \cdot (1+D^{*})$
using Good's inversion formula, Theorem~\ref{theorem_Good}.
For a similar computation; see Theorem~\ref{theorem_1_d+1}.

\item[(c)]
\cite[Corollary~4.9]{Muhle_Nadeau_Williams}
The M\"obius function of $\NC^{d}_{dk+1}$ is given by
\begin{align*}
\mu( \NC^{d}_{dk+1} )
& =
(-1)^{k} \cdot \frac{1}{2dk-k+1} \cdot \binom{2dk}{k} .
\end{align*}
Set $a^{*}(k) = s = t = 1$ and $a(k) = (-1)^{k} \cdot \frac{1}{2dk-k+1} \cdot \binom{2dk}{k}$.
Since each interval $[\hz,\pi]$ factors as the product
$\prod_{B \in \pi} \NC^{d}_{|B|}$ and that the M\"obius function factors
over Cartesian products, it is enough to prove
$\sum_{\pi \in \NC^{d}_{dk+1}} \prod_{B \in \pi} a((|B|-1)/d) = \delta_{k,0}$
where $\delta_{k,0}$ is the Kronecker delta.
The four generating functions
$A(x)$, 
$A^{*}(x)$, 
$C(x)$ and 
$C^{*}(x)$
satisfy the relations
$A(x) + x \cdot A(x)^{2d} = 1$,
$A^{*}(x) = 1/(1-x)$,
$C(x)^{d} - C(x)^{d-1} = x$, $C(0) = 1$
and $C^{*}(x) = 1/C(x)$,
which yields the result.
\end{itemize}
}
\label{example_Big_example}
\end{example}

By considering the generating function $C(x)$ in
Example~\ref{example_Big_example} parts~(a) and~(b),
we obtain the following results.
\begin{corollary}
\begin{itemize}
\item[(a)]
The number of noncrossing partitions in $\NC^{d}_{dk+1}$
that contain the singleton block $\{1\}$ is given by
\begin{align*}
\frac{1}{dk+1} \cdot \binom{dk+k}{k} .
\end{align*}
\item[(b)]
Furthermore, the number of partitions of rank $j=k-i$
in the noncrossing partition poset~$\NC^{d}_{dk+1}$
that contain $\{1\}$ as a singleton block
is given by
\begin{align*}
\frac{1}{dj+1}
\cdot
\binomial{i+dj}{i}
\cdot
\binomial{di+j-1}{j} .
\end{align*}
\end{itemize}
\end{corollary}

\subsection{Block sizes $1$ and $d+1$}
\label{subsection_Cardinality}

\begin{theorem}
\label{theorem_smallblocks}
The number of noncrossing partitions on $dk+1$ elements where
each block size is either $1$ or $d+1$
and
each block size of the dual partition is also either $1$ or $d+1$
is given by
\begin{align*}
\frac{2}{dk+2} \cdot \binom{dk+2}{k} .
\end{align*}
If we further assume that the partition contains
the singleton block $\{1\}$, the number reduces to
\begin{align*}
\frac{1}{dk+1} \cdot \binom{dk+1}{k} .
\end{align*}
\end{theorem}
\begin{proof}
In this case we use $A(x) = A^{*}(x) = 1 + x$ and $s = t = 1$.
By symmetry we have $C(x) = C^{*}(x)$
and this generating function satisfies
$C(x) = 1 + x \cdot C(x)^{d}$.
Again equation~\eqref{equation_coefficient_of_the_generating_function}
yields these two results.
\end{proof}

\begin{theorem}
Let $i+j=k$.
The number of noncrossing partitions $\pi$ in $\NC^{d}_{dk+1}$ of rank $i$ such that
both $\pi$ and its dual partition $\pi'$
only has blocks of sizes $1$ and $d+1$ is given by
\begin{align*}
\frac{dk+1}{(dj+1) \cdot (di+1)}
\cdot
\binomial{dj+1}{i}
\cdot
\binomial{di+1}{j} .
\end{align*}
\label{theorem_1_d+1}
\end{theorem}
\begin{proof}
This proof is similar to the proof of 
Example~\ref{example_Big_example}(b).
Let $D(s,t,x) = C(s,t,x) - 1$ and $D^{*}(s,t,x) = C^{*}(s,t,x) - 1$.
Rewrite the two equations in~\eqref{equation_C_and_C*}
recalling that in this case we have $A(x) = A^{*}(x) = 1+x$, as
\begin{align*}
D
& =
sx \cdot (1+D^{*})^{d} , &
D^{*}
& =
tx \cdot (1+D)^{d} .
\end{align*}
Set $x=1$ and apply Good's inversion formula using
$g_{1}(s,t) = (1+t)^{d}$
and
$g_{2}(s,t) = (1+s)^{d}$.
Let $M$ be the matrix given by
\begin{align*}
M_{1,1}
& =
g_{1}(s,t)^{i} - s \cdot \frac{\partial g_{1}}{\partial s} \cdot g_{1}(s,t)^{i-1}
= 
(1+t)^{di} , 
&
M_{1,2}
& =
- d \cdot t \cdot (1+t)^{di-1} , \\
M_{2,1}
& =
- s \cdot \frac{\partial g_{2}}{\partial s} \cdot g_{2}(s,t)^{j-1}
=
- d \cdot s \cdot (1+s)^{dj-1} , 
&
M_{2,2}
& =
(1+s)^{dj} .
\end{align*}
Its determinant is given by
\begin{align*}
\det(M)
& =
(1+s)^{dj} \cdot (1+t)^{di}
-
 d^{2} \cdot st
\cdot
(1+s)^{dj-1} \cdot (1+t)^{di-1}  .
\end{align*}
Good's inversion formula now implies
\begin{align}
\nonumber
[s^{i} t^{j}] D
& =
[s^{i-1} t^{j}] \det(M) 
=
\binomial{dj}{i-1} \cdot \binomial{di}{j}
-
d^{2} \cdot \binomial{dj-1}{i-2} \cdot \binomial{di-1}{j-1} 
\nonumber
\end{align}
and
\begin{align}
\nonumber
[s^{i} t^{j}] D \cdot D^{*}
& =
[s^{i-1} t^{j-1}] \det(M)
=
\binomial{dj}{i-1} \cdot \binomial{di}{j-1}
-
d^{2} \cdot \binomial{dj-1}{i-2} \cdot \binomial{di-1}{j-2} .
\nonumber
\end{align}
We have that
\begin{align*}
[s^{i} t^{j}] B(s,t)
& =
[s^{i} t^{j}] D \cdot D^{*}
+
[s^{i} t^{j}] D
+
[s^{i} t^{j}] D^{*}
+
[s^{i} t^{j}] 1.
\end{align*}
Hence the result follows by adding the two previous identities and
the symmetric version for $[s^{i} t^{j}] D^{*}$.
\end{proof}

\section{A tree representation behind the enumerative results}
\label{section_A_tree_representation_of_noncrossing_partitions}

In this section we introduce a two-colored labeled plane tree 
representation of a noncrossing partition which simultaneously
represents its blocks and the blocks of its
Kreweras dual. This representation helps visualize the proofs of our
enumerative results in in Section~\ref{section_Enumerative}. We consider
a plane tree to be not only a {\em topological tree}, that is, a tree
with a cyclic order on the edges around each vertex, but a topological
tree enriched with a vertex marked as root. The labeled topological tree
we use is the same as the one introduced by M\"uhle, Nadeau and
Williams~\cite[Figure~3]{Muhle_Nadeau_Williams} but the way we mark and
use a root vertex is different from the way they mark and use their root
edge to prove their enumerative results. The difference between the two
uses of the same topological trees in recursive enumeration will be
outlined at the end of this section

The idea to use labeled two-colored rooted trees to visually 
represent a factorization $(1,2,\ldots,n)=\alpha\beta$ of a cyclic
permutation into two factors first appeared in the work of Goulden and
Jackson~\cite[Theorem~2.1]{Goulden_Jackson}. The special case considered
in~\cite[Figure~3]{Muhle_Nadeau_Williams} is when $\alpha$ represents a
noncrossing partition and $\beta$ its Kreweras dual.

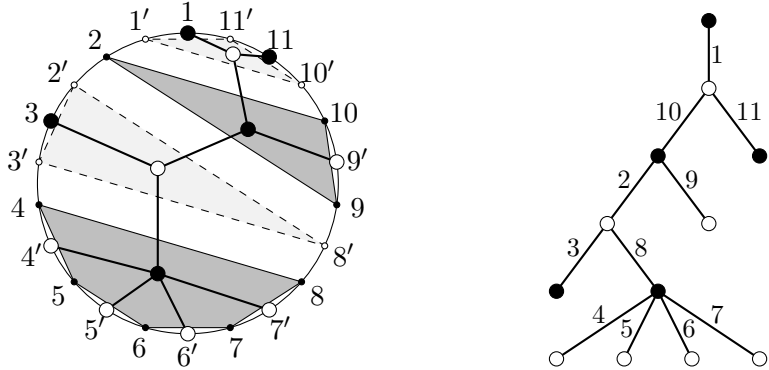
\begin{figure}[ht]
\newcommand{\app}{360/11}
\newcommand{\bpp}{90-360/11}
\newcommand{\gpp}[1]{({5*cos(\app*#1+\bpp)},{5*sin(\app*#1+\bpp)})}
\newcommand{\hpp}[1]{({5.7*cos(\app*#1+\bpp)},{5.7*sin(\app*#1+\bpp)})}
\newcommand{\bqq}{90-360/11+360/22}
\newcommand{\gqq}[1]{({5*cos(\app*#1+\bqq)},{5*sin(\app*#1+\bqq)})}
\newcommand{\hqq}[1]{({5.7*cos(\app*#1+\bqq)},{5.7*sin(\app*#1+\bqq)})}
\newcommand{\rrr}{0.24}
\begin{center}
\begin{tikzpicture}[scale = 0.40,inner sep=0.1mm]
\draw[-] (0,0) circle (5);

\draw[fill=gray!50] \gpp{2} -- \gpp{9} -- \gpp{10} -- cycle;
\draw[fill=gray!50] \gpp{4} -- \gpp{5} -- \gpp{6} -- \gpp{7} -- \gpp{8} -- cycle;
\draw[dashed,fill=gray!10] \gqq{1} -- \gqq{10} -- \gqq{11} -- cycle;
\draw[dashed,fill=gray!10] \gqq{2} -- \gqq{3} -- \gqq{8} -- cycle;

\node (blockone) at (2,1.8) {};
\node (blocktwo) at (-1,-3) {};
\node (dualblockone) at (1.5,4.3) {};
\node (dualblocktwo) at (-1,0.5) {};

\draw[-,thick] \gpp{1} -- (dualblockone);
\draw[-,thick] (dualblockone) -- \gpp{11};
\draw[-,thick] (dualblockone) -- (blockone);
\draw[-,thick] (blockone) -- \gqq{9};
\draw[-,thick] (blockone) -- (dualblocktwo);
\draw[-,thick] (dualblocktwo) -- \gpp{3};
\draw[-,thick] (dualblocktwo) -- (blocktwo);
\draw[-,thick] (blocktwo) -- \gqq{4};
\draw[-,thick] (blocktwo) -- \gqq{5};
\draw[-,thick] (blocktwo) -- \gqq{6};
\draw[-,thick] (blocktwo) -- \gqq{7};

\foreach \k in {1, 2, ..., 11}
{
\node at \hpp{\k} {$\k$};
\node at \hqq{\k} {$\k^{\prime}$};
\draw[-,fill] \gpp{\k} circle (0.1);
\draw[-,fill=white] \gqq{\k} circle (0.1);
};

\draw[-,fill] (blockone) circle (\rrr);
\draw[-,fill] (blocktwo) circle (\rrr);
\draw[-,fill=white] (dualblockone) circle (\rrr);
\draw[-,fill=white] (dualblocktwo) circle (\rrr);

\draw[-,fill] \gpp{1} circle (\rrr);
\draw[-,fill] \gpp{3} circle (\rrr);
\draw[-,fill] \gpp{11} circle (\rrr);
\draw[-,fill=white] \gqq{4} circle (\rrr);
\draw[-,fill=white] \gqq{5} circle (\rrr);
\draw[-,fill=white] \gqq{6} circle (\rrr);
\draw[-,fill=white] \gqq{7} circle (\rrr);
\draw[-,fill=white] \gqq{9} circle (\rrr);
\end{tikzpicture}
\hspace*{20 mm}
\renewcommand{\rrr}{0.21}
\begin{tikzpicture}[scale = 0.45,inner sep=0.1mm]
\node (nodeone) at (0,5) {};
\node (nodetwo) at (0,3) {};
\node (nodethree) at (-1.5,1) {};
\node (nodefour) at (1.5,1) {};
\node (nodefive) at (-3.0,-1) {};
\node (nodesix) at (0,-1) {};
\node (nodeseven) at (-4.5,-3) {};
\node (nodeeight) at (-1.5,-3) {};
\node (nodenine) at (-4.5,-5) {};
\node (nodeten) at (-2.5,-5) {};
\node (nodeeleven) at (-0.5,-5) {};
\node (nodetwelve) at (1.5,-5) {};

\draw[-,thick] (nodeone) --node[right]{\small $1$} (nodetwo);
\draw[-,thick] (nodetwo) --node[anchor=south east]{\small $10$} (nodethree);
\draw[-,thick] (nodetwo) --node[anchor=south west]{\small $11$} (nodefour);
\draw[-,thick] (nodethree) --node[anchor=south east]{\small $2$} (nodefive);
\draw[-,thick] (nodethree) --node[anchor=south west]{\small $9$} (nodesix);
\draw[-,thick] (nodefive) --node[anchor=south east]{\small $3$} (nodeseven);
\draw[-,thick] (nodefive) --node[anchor=south west]{\small $8$} (nodeeight);
\draw[-,thick] (nodeeight) --node[anchor=south east]{\small $4$} (nodenine);
\draw[-,thick] (nodeeight) --node[anchor=15]{\small $5\:$} (nodeten);
\draw[-,thick] (nodeeight) --node[anchor=165]{\small $\:6$}  (nodeeleven);
\draw[-,thick] (nodeeight) --node[anchor=south west]{\small $7$}  (nodetwelve);

\draw[-,fill] (nodeone) circle (\rrr);
\draw[-,fill=white] (nodetwo) circle (\rrr);
\draw[-,fill] (nodethree) circle (\rrr);
\draw[-,fill] (nodefour) circle (\rrr);
\draw[-,fill=white] (nodefive) circle (\rrr);
\draw[-,fill=white] (nodesix) circle (\rrr);
\draw[-,fill] (nodeseven) circle (\rrr);
\draw[-,fill] (nodeeight) circle (\rrr);
\draw[-,fill=white] (nodenine) circle (\rrr);
\draw[-,fill=white] (nodeten) circle (\rrr);
\draw[-,fill=white] (nodeeleven) circle (\rrr);
\draw[-,fill=white] (nodetwelve) circle (\rrr);
\end{tikzpicture}
\end{center}
\caption{A noncrossing partition and its plane tree representation}
\label{figure_nctree}
\end{figure}

First we outline how the same topological tree representation also arises
directly by inspecting the circle representation of the noncrossing
partition. Our more geometric approach to rooted bicolored trees representing
noncrossing partitions is a variant of the 
``dual tree representation'' discussed by Kortchemski and
Marzouk~\cite[Section 2]{Kortchemski_Marzouk} relying on a bijection
introduced by Janson and Stef\'ansson~\cite{Janson_Stefansson}. 
Consider the circle representation of a noncrossing partition $\pi$ of
$\{1,2,\ldots,n\}$. On the left hand side 
of Figure~\ref{figure_nctree} we see the circle representation of the
noncrossing partition $\pi=1|2,9,10|3|4,5,6,7,8|11$ that was introduced in
Figure~\ref{figure_example_noncrossing_partition}. The non-singleton
blocks of its Kreweras dual
partition $\pi'=1^{\prime},10^{\prime},11^{\prime}|2^{\prime},3^{\prime},8^{\prime}|4^{\prime}|5^{\prime}|6^{\prime}|7^{\prime}|9^{\prime}$
are polygons with dashed line boundaries.

Recall that a block $B$ of noncrossing
partition $\pi$ and a block $C^{\prime}$ in the dual partition $\pi^{\prime}$
are adjacent if the block $B$ contains the two elements $i$ and $j$
and the dual block $C^{\prime}$ contains the two elements $i^{\prime}$
and $(j-1)^{\prime}$.
Here $i=j$ and $i^{\prime}=(j-1)^{\prime}$ are possible.
Furthermore,
in Subsection~\ref{subsection_The_noncrossing_partition_lattice}
we defined
$\gamma(B,C^{\prime})$
to be the element in
the most in the clockwise orientation of $i$ and $j$ from the block~$B$'s
perspective.
Observe that this is the element $i$
since it has $i$ in the block $B$ and
$i^{\prime}$ in the dual block $C^{\prime}$.

\begin{definition}
\label{definition_toptree}
The {\em labeled topological tree representing the noncrossing partition
$\pi\in\NC_{n}$} is a two vertex-colored tree with $n$ labeled
edges, defined as follows:
\begin{enumerate}  
\item
The black nodes represent the blocks of $\pi$, the white
nodes represent the blocks of $\pi'$.
\item
If a block $B$ of $\pi$ and a block $C^{\prime}$ in the dual partition
are adjacent then the associated nodes are connected with an edge.
\item
Label the edge between the two nodes representing the two adjacent blocks
$B$ and $C^{\prime}$
with the canonical element $\gamma(B,C^{\prime})$.
\end{enumerate}  
\end{definition}    
The labeled topological tree associated to the noncrossing partition $\pi$ is
shown on the left hand side of Figure~\ref{figure_nctree} in such a way
that the vertex representing a non-singleton (dual) block is inside the polygon
representing the same block, whereas the vertices representing the
singleton (dual) blocks appear at correct place representing the
respective point $i$ ($i'$) on the circle. So far our representation is
essentially the same as the dual tree representation
in~\cite[Section~2]{Kortchemski_Marzouk}, the only difference being that they draw the
edge labels near the respective black vertex. 
We place the same labels in the middle of each vertex, and we think of the
label $i$ as the representative of $i$ in the block of $\pi$ associated
to the black end of the edge and as the representative of $i'$ in the
block of $\pi'$ associated to the white end of the same edge.
\begin{remark}
{\rm The careful reader should notice that the labeled topological trees
  shown in~\cite[Figure~2]{Kortchemski_Marzouk} and
  in~\cite[Figure~3]{Muhle_Nadeau_Williams} are the
  same as ours after taking the mirror image of the figure. The topological
  graph structure in the cited sources is defined in terms of the permutations
  associated to the noncrossing partition and to its Kreweras dual.  
}    
\end{remark}  

As in~\cite[Section~2]{Kortchemski_Marzouk}, the following statement is
a direct consequence of this construction.
\begin{proposition}
\label{proposition_injective1}  
The operation associating to each noncrossing partition $\pi\in\NC_{n}$ its labeled topological tree is injective.
\end{proposition}  

Next we turn our labeled topological tree into a {\em (rooted) plane
  tree}. This is where our choices differ
from~\cite[Section~2]{Kortchemski_Marzouk} substantially, as we define
the root and its leftmost child (called sometimes ``corner''
in~\cite[Section~2]{Kortchemski_Marzouk}) in a different
fashion. The way of turning a topological tree into a plane tree is
implicit in~\cite{Muhle_Nadeau_Williams}, but also different. 

\begin{definition}
\label{definition_planetree}  
Given a noncrossing partition $\pi\in\NC_{n}$ we define its {\em labeled
  plane tree} as the drawing of its labeled topological tree in the
plane subject to the following rules.  
\begin{enumerate}
\item The root of the tree represents the block $B_{1}$ of $\pi$ containing the element~$1$.
\item The leftmost child of each vertex is connected to it by the edge whose
  label is the least among the labels of all edges connecting the vertex
  to its children.  
\end{enumerate}  
\end{definition} 
The labeled plane tree representation of the noncrossing partition shown
in Figure~\ref{figure_nctree} is on the right hand side of the figure.
Note that we may omit the mention of the
two-coloring of the vertices, because assigning the color black to the root
determines the two-coloring uniquely.
As a consequence of Proposition~\ref{proposition_injective1}
we obtain the following corollary.
\begin{corollary}
  \label{corollary_injective2}
The operation associating to each noncrossing partition $\pi\in\NC_{n}$
its labeled plane tree is injective.
\end{corollary}  

Next we show that keeping track of the labels is not necessary either,
as it may be uniquely reconstructed from the plane tree structure. The
following statement plays a key role in seeing this. Since the vertices
represent (dual) blocks of a noncrossing partition, by a slight abuse of
notation we will denote the label on the edge connecting the vertices
$u$ and $v$ also by $\gamma(u,v)$.
\begin{proposition}
\label{proposition_labelkey}  
The labeling of a labeled plane tree representing a noncrossing
partition $\pi\in\NC_{n}$ has the following properties:
\begin{enumerate}
\item The labels on the edges are the elements of the set
  $\{1,2,\ldots,n\}$ each appearing exactly once.
\item For each vertex in the tree, the labels on the edges connecting it
  to its children increase left to right.
\item For each nonroot black (white) vertex $u$ in the tree the largest
  (least) label appearing on an edge connecting $u$ to another vertex is
  on the edge connecting $u$ to its parent in the tree.
\item For each vertex in the tree the set of all labels appearing in its
  subtree form an interval $[i,j]=\{i,i+1,\ldots,j\}$ of consecutive
  integers.
\item For each black (white) vertex in the tree, the least (largest)
  label in its subtree labels the edge connecting the root of the
  subtree to its leftmost (rightmost) child.   
\item If $v_{1}$ and $v_{2}$ are children of the same vertex $v$ and $v_{1}$
  is to the left from $v_{2}$ then $\gamma(v,v_{1})$ is smaller than the label
  on any edge in the subtree of $v_{2}$, and $\gamma(v,v_{2})$ is greater
  than the label on any edge in the subtree of $v_{1}$.  
\end{enumerate}  
\end{proposition}  
\begin{proof}
In our proof we use the fact that each black vertex represents a block  
$\{i_{1},i_{2},\ldots,i_{r}\}$ of a noncrossing partition, here we assume
$i_{1}<i_{2}<\cdots<i_{r}$. Similarly, each white vertex represents a dual
block $\{j_{1}',j_{2}',\ldots,j_{s}'\}$, here we assume
$j_{1}<j_{2}<\cdots<j_{s}$. 

Property~(1) is a consequence of
Definition~\ref{definition_toptree}. It is also a consequence of this
definition that for each vertex in the tree, the labels on the edges
connecting it to its children follow the counterclockwise order. By
part~(2) of Definition~\ref{definition_planetree} each such cyclic list of
labels begins with its least element, this implies property~(2).
Next we prove part~(3). There is nothing to show for the root
vertex, which represents a block $\{i_{1},i_{2},\ldots,i_{r}\}$ satisfying $i_{1}=1$. 
Any other black vertex represents a block $\{i_{1},i_{2},\ldots,i_{r}\}$ 
where $1<i_{1}<\cdots<i_{r}$. The edge connecting this vertex to its
parent crosses the side $\{i_{1},i_{r}\}$ because this line separates the
polygon from the point $1$. The label of this edge is $i_{r}$. The proof
of property~(3) for white vertices is similar. (Note that a
white vertex can never be the root of the entire tree). A white vertex
represents a dual block $\{j_{1}',j_{2}',\ldots,j_{s}'\}$, where 
$j_{1}<j_{2}<\cdots<j_{s}$. Once again the edge
connecting this vertex to its parent crosses the side $\{j_{1}',i_{s}'\}$
because this line separates the polygon from the point $1$.

We prove the remaining parts~(4), (5) and (6) at once. Consider first a
nonroot black vertex $u$ representing the block  $\{i_{1},i_{2},\ldots,i_{r}\}$ where 
$1<i_{1}<i_{2}<\cdots<i_{r}$. The set of labels appearing in its subtree is 
$\{i_{1},i_{1}+1,\ldots,i_{r}-1\}$, because the points appearing on the
counterclockwise (half open) arc $[i_{1},i_{r})$ in our circle representation are
$i_{1},i_{1}',(i_{1}+1),(i_{1}+1)',\ldots,i_{r}-1, (i_{r}-1)'$. (Keep in mind that
$i_{r}$ is the label of the edge connecting our vertex to its parent.)
  Furthermore, the set $\{i_{1},i_{2},\ldots,i_{r-1}\}$ is 
used to label the edges connecting $u$ to its children. The labels
strictly between $i_{t}$ and $i_{t+1}$ (for $1\leq t\leq r-1$) are the
ones appearing in the subtree of the child $v$ of $u$ satisfying
$\gamma(u,v)=i_{t}$.  

For the root vertex representing the block  $\{i_{1},i_{2},\ldots,i_{r}\}$ where 
$1=i_{1}<i_{2}<\cdots<i_{r}$, we need to amend slightly the above
reasoning. The subtree of the root is the entire tree, containing all
labels. The set $\{i_{1},i_{2},\ldots,i_{r}\}$ is 
used to label the edges connecting the root to its children. The labels
strictly between $i_{t}$ and $i_{t+1}$ (for $1\leq t\leq r-1$) are the
ones appearing in the subtree of the child $v$ of $u$ satisfying
$\gamma(u,v)=i_{t}$, and the labels appearing in the subtree of the last
child of the root are the larger than all other labels. 

For white vertices parts~(4), (5) and (6) may be shown in a similar
fashion. A white vertex represents a dual block
$\{j_{1}',j_{2}',\ldots,j_{s}'\}$ 
where $j_{1}<j_{2}<\cdots<j_{s}$ holds. The label of
the edge connecting our white vertex to its parent is $j_{1}'$, and the
labels appearing in the subtree of our vertex are
$\{j_{1}+1,\ldots,j_{s}\}$, because the points appearing on the
counterclockwise (half open) arc $(j_{1}',j_{s}']$ in our circle representation are
$j_{1}+1,(j_{1}+1)',(j_{1}+2),(j_{1}+2)',\ldots,j_{s},j_{s}'$. The set
$\{j_{2},j_{3}\ldots,j_{s}\}$ is used to label the edges connecting~$u$ to its
children. The labels 
strictly between $j_{t-1}$ and $j_{t}$ for $2\leq t\leq s$ are the ones
appearing in the subtree of of the child $v$ of $u$ satisfying
$\gamma(u,v)=j_{s}$.    
\end{proof}  
\begin{definition}
  We call the labeled plane trees satisfying the conditions stated in
  Proposition~\ref{proposition_labelkey}  above {\em properly labeled},
  and we denote the set of properly labeled plane trees on $n+1$
  vertices by $\ptree_{n}$.  
\end{definition}  

\begin{proposition}
\label{proposition_unique}  
Given a plane tree with $n+1$ vertices, there is a unique way to
properly label its edges.
\end{proposition}
\begin{proof}
We prove by induction on the number of labels the following statement:
given a vertex, and the set of labels on the edges of its subtree, there
is at most one way to place the labels on the edges. The statement is
obviously true for the empty set, and we will prove the induction step
for the subtree of a black vertex $u$: the case of the subtree of a
white vertex is completely analogous.

Assume that the set of children of $u$ is $\{v_{1},\ldots,v_{k}\}$, listed in
the left to right order, and the label of the edge $\{v,v_{t}\}$ is $i_{t}$
for $t=1,2,\ldots,k$. By property~(2) we have $i_{1}<i_{2}<\cdots <i_{k}$ and
by property~(5) the label $i_{1}$ is the least element of the set of all
labels appearing in the subtree. By property (4) the set of all labels
in the subtree is a set of consecutive integers, the same holds for the
subtrees of $v_{1}$, $v_{2}$, \ldots and $v_{k}$. Combining properties~(3)
through (5) we obtain that the labels in the subtree of $u$ must be
listed in the following increasing order: $i_{1}$, labels in the subtree
of $v_{1}$, $i_{2}$, labels in the subtree of $v_{2}$, \ldots, $i_{k}$, labels
in the subtree of $v_{k}$. (Keep in mind that each $v_{t}$ is a white vertex.) 
There is exactly one way to partition the set of labels and assign them
to the subtrees of each $v_{t}$ and to the edges connecting $u$ to its
children. 
\end{proof}  
\begin{theorem}
\label{theorem_treerep}
The operation assigning to each noncrossing partition $\pi\in\NC_{n}$
its properly labeled plane tree $T\in\ptree_{n}$ is a bijection.
\end{theorem}  
\begin{proof}
The operation is an injection by Proposition~\ref{proposition_injective1} and
Corollary~\ref{corollary_injective2}. By Proposition~\ref{proposition_unique}, no
information is lost either if we remove the labels of the edges. The
surjectivity now follows from the well known fact that the 
noncrossing partitions of $\{1,2, \ldots, n\}$ and the plane
trees on $n+1$ vertices are enumerated by the same Catalan number.  
\end{proof}

Since the degree of each vertex in our tree representations equals the
number of elements of the block it represents, we obtain the following
corollaries to
Theorem~\ref{theorem_treerep}. Following~\cite[Remark~4.2]{Muhle_Nadeau_Williams},we
call a rooted plane tree {\em $d$-divisible} if the number of children
of each vertex is divisible by $d$.  

\begin{corollary}
The cardinality of the poset $\NC^{d}_{n}$ is
the number of rooted plane trees on $n+1$ vertices
such that the degree of each vertex is congruent to $1$ modulo $d$.
\end{corollary}  
\begin{corollary}
\label{corollary_bigblocks}  
The number of noncrossing partitions in $\NC^{d}_{dk+1}$
that contain the singleton block~$\{1\}$ is the number of $d$-divisible
rooted plane trees on $dk+1$ vertices.
\end{corollary}  
Indeed, an element of $\NC^{d}_{dk+1}$ contains the singleton block
$\{1\}$ if and only if the root of the corresponding plane tree on
$dk+2$ vertices has a single child. Such plane trees correspond
bijectively to the plane trees on $dk+1$ vertices obtained by removing
the root and designating its only child as the new root. Using the same
idea we also obtain the following corollary.
\begin{corollary}
\label{corollary_smallblocks}    
The number of noncrossing partitions in $\NC^{d}_{dk+1}$ having only
blocks and dual blocks of sizes $1$ and $d+1$ and 
containing the singleton block $\{1\}$ is the number of rooted $d$-ary
trees on $dk+1$ vertices. 
\end{corollary}
Hence the second statement in Theorem~\ref{theorem_smallblocks} is a
direct consequence of the well-known fact that the number of rooted
$d$-ary trees is a Fuss--Catalan number. We may also derive the second
statement of 
Example~\ref{example_Big_example}(a)
from
Corollary~\ref{corollary_bigblocks}  by observing that the number of
$d$-divisible rooted plane trees on $dk+1$ vertices is the same as the
number of $(d+1)$-ary trees with $d+1$ nonleaf vertices. A bijective
proof of this statement is outlined
in~\cite[Remark~4.2]{Muhle_Nadeau_Williams}.

\begin{remark}
{\rm
Proposition~6.2.1 in~\cite{Stanley_EC2}
presents six classes of objects that are in bijection with each other.
See also the historical remarks in the Notes in~\cite[Chapter~6]{Stanley_EC2}.
We add a seventh class to this collection.
Let $S$ be a subset of the positive integers $\Ppp$
and let $n$ and $m$ be two positive integers.
The set of noncrossing partitions $\pi$ in the set $\NC_{n}$ such that
for all blocks $B$ of $\pi$ and the blocks $C^{\prime}$ of~$\pi^{\prime}$
the block sizes belongs
to the set $\{1\} \cup (S+1)$,
the singleton block $\{1\}$ belongs to the partition $\pi$
and
there are $m+1$ singleton blocks total among the two 
partitions $\pi$ and the dual $\pi^{\prime}$.
For instance in the case $S = \{2,3\}$, $n=6$ and $m=4$, the noncrossing partition $1|2|3,4,5,6$
corresponds to the tree in~\cite[Proposition~6.2.1, part~(i)]{Stanley_EC2}.
}
\end{remark}

We conclude this section with visually highlighting the main difference between
the main idea used to obtain enumerative results
in~\cite{Muhle_Nadeau_Williams} and in our work. Each
labeled two-colored topological tree is decomposed into a pair of
$d$-divisible plane trees in~\cite{Muhle_Nadeau_Williams} by removing
the {\em root edge} labeled $1$. The roots of the two trees have
opposite color. Our approach corresponds to splitting the {\em root
  vertex} $v$, that is, the black vertex incident to the edge labeled
$1$ into $\deg(v)$ copies. Thus we obtain several labeled topological
trees of the same structure, whose root vertex has degree $1$. After
identifying each new root vertex with its only child, we obtain a
collection of $\deg(v)$ $d$-divisible rooted plane trees. The use of the
Good inversion formula seems unavoidable with either approach: in the
paper of M\"uhle, Nadeau and Williams~\cite{Muhle_Nadeau_Williams} this
work is done in a cited paper of Goulden and Jackson~\cite{Goulden_Jackson}.

\section{The antipode of the poset $NC^{d}_{n}$}
\label{section_antipode}

Let ${\mathcal P}$ denote the linear span of all the isomorphism classes of finite posets
with minimal element~$\hz$ and maximal element~$\ho$.
The vector space ${\mathcal P}$ forms a Hopf algebra. That is, there is
a product, a coproduct, a unit, a counit and an antipode.
For more details; see~\cite{Ehrenborg,Schmitt}.
The antipode is described by the following expression due to
Schmitt~\cite[Theorem~6.1]{Schmitt}:
\begin{align}
S(P)
& =
\sum_{k \geq 0}
\sum_{c}
(-1)^{k}
\cdot
[x_{0},x_{1}] \cdot [x_{1},x_{2}] \cdots [x_{k-1},x_{k}] ,
\label{equation_Schmitt_antipode}
\end{align}
where the sum is over all chains
$c = \{\hz = x_{0} < x_{1} < \cdots < x_{k} = \ho\}$
taking $k$ steps.
Observe that this is a generalization of Philip Hall's
formula for the M\"obius function.
The disadvantage with the above sum is that it yields
an expression that may contain a lot of cancellations.

A {\em hypergraph $G$} on a vertex set ${\mathcal V}$
is a collection ${\mathcal E}$ of subsets of ${\mathcal V}$
such that the subsets $E \in {\mathcal E}$ have cardinality at least $2$.
The subsets in ${\mathcal E}$ are called edges.
A hypergraph is connected if for every pair of vertices $x$ and $y$ there is a
sequence of edges $E_{1}, E_{2}, \ldots, E_{k} \in {\mathcal E}$ such that
$x \in E_{1}$, $y \in E_{k}$ and the intersection
$E_{i} \cap E_{i+1}$ is nonempty for all $1 \leq i \leq k-1$.
A connected hypergraph is a hypertree if 
$\sum_{E \in {\mathcal E}} (|E|-1) = |{\mathcal V}| - 1$.
Note that if all the edges have cardinality $2$ we have a tree.

Einziger~\cite{Einziger_dissertation,Einziger_preprint}
proved the following result for the antipode of the noncrossing partition lattice.
A proof using a sign reversing involution was given in~\cite{Ehrenborg_Happ}.

\begin{theorem}[Einziger]
The antipode of the noncrossing partition lattice $\NC_{n}$ is given by the sum
\begin{align*}
S(\NC_{n})
& =
\sum_{T}
(-1)^{|T|}
\cdot
\prod_{E \in T} \NC_{|E|} ,
\end{align*}
where $T$ ranges over all noncrossing hypertrees on the set $[n]$.
\end{theorem}

We now extend this result to the noncrossing partition poset $\NC^{d}_{n}$.

\begin{theorem}
Let $n$ and $d$ be positive integers such that $n \equiv 1 \bmod d$.
Then the antipode of the noncrossing partition poset $\NC^{d}_{n}$ is given by the sum
\begin{align*}
S(\NC^{d}_{n})
& =
\sum_{T}
(-1)^{|T|}
\cdot
\prod_{E \in T} \NC^{d}_{|E|} ,
\end{align*}
where $T$ ranges over all noncrossing hypertrees on the set $[n]$
such that each edge $E$ of $T$ satisfies $|E| \equiv 1 \bmod d$.
\label{theorem_antipode}
\end{theorem}

Recall that
in Subsection~\ref{subsection_The_noncrossing_partition_lattice}
for a block $B$ of a partition $\pi$ and
an adjacent block $C'$ in the dual partition $\pi'$, we defined the notion
$\gamma(B,C')$ to be the vertex of $B$ adjacent to the block $C'$
that is the most negative orientation from the block $B$'s perspective.

Next, for two noncrossing partitions~$\pi$ and~$\sigma$ in
in $\NC^{d}_{n}$ such that $\pi < \sigma$ define the
hypergraph $\varphi(\pi,\sigma)$ as follows.
For each block $C^{\prime}$ in the dual partition
and a maximal collection of blocks
$B_{1}, B_{2}, \ldots, B_{j}$ of~$\pi$, where $j \geq 2$,
which are adjacent to the block $C^{\prime}$
and are all contained in one block of the partition~$\sigma$, add
the edge
$\{\gamma(B_{1},C^{\prime}),\gamma(B_{2},C^{\prime}), \ldots, \gamma(B_{j},C^{\prime})\}$
to the hypergraph $\varphi(\pi,\sigma)$.
It is important to note that $j$ satisfies the congruence
$j \equiv 1 \bmod d$.

For a chain $c = \{\hz = \pi_{0} < \pi_{1} < \cdots < \pi_{r} = \ho\}$
in the noncrossing partition poset $\NC^{d}_{n}$
define the hypergraph $\varphi(c)$ to be the union
$\cup_{i=1}^{r} \varphi(\pi_{i-1},\pi_{i})$.

Similar to \cite[Lemma~3.3 and Corollary~3.5]{Ehrenborg_Happ}
we have the next result whose proof we omit.
\begin{lemma}
For a chain $c = \{\hz = \pi_{0} < \pi_{1} < \cdots < \pi_{r} = \ho\}$ in the poset $\NC^{d}_{n}$
the following isomorphism holds
\begin{align*}
\prod_{i=1}^{r} [\pi_{i-1},\pi_{i}] 
& \cong
\prod_{E \in \varphi(c)} \NC^{d}_{|E|} .
\end{align*}
\label{lemma_Ehrenborg_Happ_Lemma_3.3_Corollary_3.5}
\end{lemma}
Finally, their Lemma~3.6 becomes:
\begin{lemma}
For a chain $c$ in the poset $\NC^{d}_{n}$
the hypergraph $\varphi(c)$ is a noncrossing hypertree where each edge
$E$ satisfies $|E| \equiv 1 \bmod d$.
\label{lemma_Ehrenborg_Happ_Lemma_3.6}
\end{lemma}
The only difference in the proof of this lemma is the factor of $d$
in the following chain of equalities:
\begin{align*}
\sum_{E \in \varphi(c)} (|E|-1)
=
d \cdot \sum_{E \in \varphi(c)} \rho\left(\NC^{d}_{|E|}\right)
=
d \cdot \sum_{i=1}^{r} \rho([\pi_{i-1},\pi_{i}])
=
d k = n-1 .
\end{align*}
Finally, Proposition~4.3 in~\cite{Ehrenborg_Happ} states:
\begin{proposition}[Ehrenborg--Happ]
For a noncrossing hypertree $H$ on $n$ elements with $r$ edges,
the following alternating sum holds:
\begin{align*}
\sum_{c \in \varphi^{-1}(H)} (-1)^{\ell(c)}
& =
(-1)^{r} .
\end{align*}
\label{proposition_fiber}
\end{proposition}

\begin{proof}[Proof of Theorem~\ref{theorem_antipode}.]
By combining the Schmitt expression for the antipode \eqref{equation_Schmitt_antipode}
and Proposition~\ref{proposition_fiber} the result follows.
\end{proof}

By combining Example~\ref{example_Big_example}(b)
and Theorem~\ref{theorem_antipode}
using Philip Hall's formula for the M\"obius function,
we obtain:

\begin{corollary}
Let $k$ and $d$ be nonnegative integers such that $d \geq 1$.
Then the following holds
\begin{align*}
(-1)^{k} \cdot \frac{1}{2dk-k+1} \cdot \binom{2dk}{k}
& =
\sum_{T} (-1)^{|T|} ,
\end{align*}
where $T$ ranges over all noncrossing hypertrees on the set $[dk+1]$
such that each edge $E$ of $T$ satisfies $|E| \equiv 1 \bmod d$.
\end{corollary}

\section{An $EL$-labeling of the poset $NC^{d}_{n}$}
\label{section_EL-labeling_of_the_poset}

We define an edge labeling of the noncrossing partition poset $\NC^{d}_{n}$ as follows.
Let $\pi \coveredby \sigma$ be a cover relation is the poset $\NC^{d}_{n}$.
By Theorem~\ref{theorem_NC_structure}
there are $d+1$ blocks of $\pi$ that are joined to form a block of $\sigma$.
Let us denote the joined blocks by $B_{1}, B_{2},\ldots, B_{d+1}$
and assume that the inequalities $\min(B_{1}) < \min(B_{2}) < \cdots < \min(B_{d+1})$ hold.
Then let the edge label $\lambda(\pi,\sigma)$ be given by
\begin{align}
\lambda(\pi,\sigma)
& =
\max(\{i \in B_{1} : i < \min(B_{2})\}) .
\label{equation_the_edge_labeling}
\end{align}
This edge labeling generalizes the edge labeling of $\NC_{n}$ introduced
by Stanley~\cite{Stanley}. 
For a maximal chain
${\bf m} = \{\hz = \pi_{0} \coveredby \pi_{1} \coveredby \cdots  \coveredby \pi_{k} = \ho\}$
in the poset $\NC^{d}_{n}$
define the labeling of the chain ${\bf m}$ to be
\begin{align*}
\lambda({\bf m})
& =
(\lambda(\pi_{0},\pi_{1}),
\lambda(\pi_{1},\pi_{2}),
\ldots,
\lambda(\pi_{k-1},\pi_{k})) .
\end{align*}

\begin{proposition}
The edge labeling in equation~\eqref{equation_the_edge_labeling}
is the same as the edge labeling
of 
M\"uhle, Nadeau and Williams
in~\cite[Theorem~5.3]{Muhle_Nadeau_Williams}.
\end{proposition}
\begin{proof}
In the paper~\cite{Muhle_Nadeau_Williams}
a block $B = \{b_{1} < b_{2} < \cdots < b_{k}\}$
in a noncrossing partition $\pi$ is encoded by the
cyclic permutation
$(b_{1}, b_{2}, \ldots, b_{k})$.
Hence the noncrossing partition $\pi$ is encoded by a permutation~$\overline{\pi}$
on the underlying set of elements.
Furthermore, given a cover relation $\pi < \sigma$
then there is a $d+1$ cycle $\tau$ such that
$\overline{\pi} \circ \tau = \overline{\sigma}$.
Theorem~5.3 in~\cite{Muhle_Nadeau_Williams} defines the edge
label to be the smallest element of the cycle $\tau$.

Suppose you want to merge $d+1$ blocks of the noncrossing partition $\pi$
to obtain the partition $\sigma$
and these are 
$B_{1} = \{b_{1,1} < b_{1,2} < \cdots < b_{1,j_{1}}\},
\ldots,
B_{d+1} = \{b_{d+1,1} < b_{d+1,2} < \cdots < b_{d+1,j_{d+1}}\}$.
Assume furthermore that
$\min(B_{1}) < \min(B_{2}) < \cdots < \min(B_{d+1})$.
These $d+1$ blocks correspond to the $d+1$ cycles
\begin{align*}
(b_{1,1}, b_{1,2}, \ldots, b_{1,j_{1}})
(b_{2,1}, b_{2,2}, \ldots, b_{2,j_{2}})
\cdots
(b_{d+1,1},  b_{d+1,2}, \ldots, b_{d+1,j_{d+1}}) .
\end{align*}
Now we are looking for a cyclic permutation $\tau=(t_{1},t_{2},\ldots, t_{d+1})$
such that the above permutation composed by $\tau$ on the right
merges the blocks and
yields a single cycle in which all elements are listed in increasing order.
Let us determine $t_{1}$, the least element of this cycle.
This element can be found just by looking at the first two blocks.
Indeed, there is an index $j \leq i_{1}$ such that 
$b_{1,1}, b_{1,2}, \ldots b_{1,j}$ are the $j$ least elements of the union of the $d+1$ blocks,
but the next least element is $b_{2,1}=\min(B_2)$. In that case 
we must have $t_{1} = b_{1,j}$ and $t_{2}=b_{2,i_{2}}$ (the latter is the largest element of $B_2$).
This agrees with the labeling in~\eqref{equation_the_edge_labeling}.
\end{proof}

Define a {\em $d$-parking function} to be a list 
$(a_{1}, a_{2}, \ldots, a_{k})$
of positive integers such that
when the list is ordered
$(a_{(1)} \leq a_{(2)} \leq \cdots \leq a_{(k)})$,
it satisfies the inequality
$a_{(i)} \leq d \cdot (i-1) + 1$.
This definition is a shift from the definitions occurring in the two papers~\cite{Stanley_k-parking,Yan},
where the entries are nonnegative integers.

\begin{theorem}[M\"uhle--Nadeau--Williams]
The map of sending a maximal chain ${\bf m}$ of $\NC^{d}_{dk+1}$
to its list of labels $\lambda({\bf m})$
is a bijection between all 
maximal chains of $\NC^{d}_{dk+1}$
and all $d$-parking functions of length $k$.
\label{theorem_list_of_labels_bijection}
\end{theorem}

We also remark that a similar bijection appears in the work of
Irving and Rattan~\cite[Proposition~1.6]{Irving_Rattan}, generalizing
the result of Biane~\cite[Theorem~3.1]{Biane}.
Next we show that the order complex of the poset $\NC^{d}_{n}$
is shellable by producing an $EL$-labeling;
see~\cite[Definition 3.2.1]{Wachs_Poset_topology}.

\begin{theorem}
The labeling
$\lambda^{*}(\pi,\sigma) = |\pi| - \lambda(\pi,\sigma)$
is an $EL$-labeling of the noncrossing partition poset~$\NC^{d}_{n}$.
\label{theorem_EL}
\end{theorem}
\begin{proof}
We first show that the poset $\NC^{d}_{n}$ has a unique rising
chain. Consider a rising maximal chain
\begin{align*}
  {\bf m} & = \{\hz = \pi_{0} \coveredby \pi_{1} \coveredby \cdots \coveredby \pi_{k} = \ho\} .
\end{align*}
Since ${\bf m}$ is a maximal chain we have
$|\pi_{i}|=|\pi_{i+1}|+d$ for $0 \leq i \leq k-1$
using Theorem~\ref{theorem_NC_structure}.
Hence the rising condition 
$\lambda^{*}(\pi_{i-1},\pi_{i}) \leq \lambda^{*}(\pi_{i},\pi_{i+1})$
implies
$\lambda(\pi_{i-1},\pi_{i}) \geq \lambda(\pi_{i},\pi_{i+1}) + d$.
As a consequence the $d$-parking function
$\vec{a} = (a_{1}, a_{2}, \ldots, a_{k})$ associated to ${\bf m}$ must satisfy
\begin{align*}
(a_{1}, a_{2}, \ldots, a_{k}) & \geq ((k-1)d+1, (k-2)d+1, \ldots,d+1, 1)
\end{align*}
coordinate-wise. Rearranging the entries of $\vec{a}$ into increasing
order, the resulting vector satisfies
\begin{align*}
(a_{(1)}, a_{(2)}, \ldots, a_{(k)}) & \geq (1, d+1,\ldots, (k-1)d+1)
\end{align*}
coordinate-wise.
By the definition of a $d$-parking function
the above lower bound for $\vec{a}$ is also an upper bound.
That is, equality holds.
Hence the only rising chain ${\bf m}$ is the one
associated to the $d$-parking function
$((k-1)d+1, (k-2)d+1, \ldots,d+1, 1)$.

Next we show that the unique rising chain of $\NC^{d}_{n}$
described above is also lexicographically first. Consider any other
maximal chain 
${\bf n} = \{\hz = \sigma_{0} \coveredby \sigma_{1} \coveredby \cdots \coveredby \sigma_{k} = \ho\}$
associated to the $d$-parking function $\vec{b} = (b_{1}, b_{2}, \ldots,
b_{k})$. Let $i$ be the least index such that $a_{i}\neq b_{i}$
holds. As a consequence of $a_{1}=b_{1}$, $a_{2}=b_{2}$, \ldots,
$a_{i-1}=b_{i-1}$, the first $i-1$ coordinates of $\vec{b}$ are also its
$i-1$ largest coordinates, and $b_{i}$ must satisfy $b_{i}\leq
(k-i+1)d$. Hence
$\lambda^{*}(\sigma_{i},\sigma_{i+1})>\lambda^{*}(\pi_{i},\pi_{i+1})$ must
be satisfied, while
$\lambda^{*}(\sigma_{j},\sigma_{j+1})=\lambda^{*}(\pi_{j},\pi_{j+1})$ holds
for all $j<i$. Therefore ${\bf m}$ precedes ${\bf n}$ in the
lexicographic order.

As noted by Stanley~\cite{Stanley} about the $d=1$ case, the statement
is easily generalized to an arbitrary interval of $\NC^{d}_{n}$ using
the fact that any such interval is a direct product of smaller copies of
$\NC^{d}_{n_{i}}$
\end{proof}

\begin{corollary}
The number of $d$-parking functions $(a_{1}, a_{2}, \ldots, a_{k})$
such that $a_{i} \leq a_{i+1} + d-1$ for all indices $i$ is given by
$\binom{2dk}{k}/(2dk-k+1)$.
\label{corollary_falling_chains}
\end{corollary}
\begin{proof}
This result follows from the formula for the M\"obius function of the
poset $\NC^{d}_{dk+1}$ stated in Example~\ref{example_Big_example}(b)
and Theorem~\ref{theorem_EL}.
Recall that the M\"obius function times the sign $(-1)^{k}$ enumerates
the number of falling chains.
\end{proof}

\section{Parking trees representing maximal chains in $NC^{d}_{n}$}
\label{section_tree_representation}

Parsing a maximal chain of $\NC^{d}_{dk+1}$ amounts to describing a process
of repeatedly merging $d+1$ blocks of a partition. In this
section we describe how parking trees provide a visual record
of this process.

The first bijection between labeled trees and parking functions was
introduced by Kreweras~\cite{Kreweras_suites}. A simpler construction was
presented by Knuth~\cite[Section~6.4, Exercise 31, solution (c)]{Knuth}, a
detailed description of his bijection may be found in the book chapter written
by Yan~\cite[Section~13.2.2]{Yan_chapter}. Here we follow a {\em different
approach} from another section of Yan's
work~\cite[Section~13.2.3]{Yan_chapter}, which can be easily generalized
to $d$-parking functions. For a key difference between the two
approaches see Remark~\ref{remark_ourtrees} below.

\begin{definition}
  \label{definition_parking_tree}
A {\em $d$-parking tree} is a labeled rooted plane tree on $n=dk+1$ vertices, such that
the number of children of each vertex is a multiple of $d$ and the
labeling satisfies the following conditions:
\begin{enumerate}
\item The label of the root is $\infty$.
\item The label of any other vertex is of the form $i_{j}$ where $i\in
  \{1,2,\ldots,k\}$ and $j\in \{1,2,\ldots,d\}$.
\item For each fixed $i$, the vertices labeled $i_{1},\ldots,i_{d}$ are
  consecutive labels in the left-to-right order of the same parent.  
\item Each $i\in \{1,2,\ldots,k\}$ is only used to label one $d$-element
  set of siblings.  
\item If $i_{j}$ and $i'_{j'}$ are children of the same vertex and $i<i'$
  holds then the vertex labeled $i_{j}$ is to the {\em right} of the vertex
  labeled $i'_{j'}$. 
\end{enumerate}  
\end{definition}

\begin{figure}[ht]
\newcommand{\shift}{5.5}
\begin{center}
\begin{tikzpicture}[scale = 1.3, inner sep=0.6mm]
\node (a) at (0,5) {$\infty$};
\node (b) at (-1.5,4) {$4_{1}$};
\node (c) at (-0.5,4) {$4_{2}$};
\node (d) at (0.5,4) {$2_{1}$};
\node (e) at (1.5,4) {$2_{2}$};
\node (f) at (-1,3) {$1_{1}$};
\node (g) at (0,3) {$1_{2}$};
\node (j) at (1,3) {$3_{1}$};
\node (k) at (2,3) {$3_{2}$};
\node (h) at (-1.5,2) {$5_{1}$};
\node (i) at (-0.5,2) {$5_{2}$};

\draw[-,line width=1pt] (a) -- (b);
\draw[-,line width=1pt] (a) -- (c);
\draw[-,line width=2pt] (a) -- (d);
\draw[-,line width=2pt] (a) -- (e);
\draw[-,line width=2pt] (c) -- (f);
\draw[-,line width=2pt] (c) -- (g);
\draw[-,line width=1pt] (f) -- (h);
\draw[-,line width=1pt] (f) -- (i);
\draw[-,line width=2pt] (e) -- (j);
\draw[-,line width=2pt] (e) -- (k);
\end{tikzpicture}
\hspace*{15mm}
\begin{tikzpicture}[scale = 1.3, inner sep=0.6mm]
\node (a) at (0,5) {$1$};
\node (b) at (-1.5,4) {$2$};
\node (c) at (-0.5,4) {$3$};
\node (d) at (0.5,4) {$8$};
\node (e) at (1.5,4) {$9$};
\node (f) at (-1,3) {$4$};
\node (g) at (0,3) {$7$};
\node (j) at (1,3) {$10$};
\node (k) at (2,3) {$11$};
\node (h) at (-1.5,2) {$5$};
\node (i) at (-0.5,2) {$6$};

\draw[-,line width=1pt] (a) -- (b);
\draw[-,line width=1pt] (a) -- (c);
\draw[-,line width=2pt] (a) -- (d);
\draw[-,line width=2pt] (a) -- (e);
\draw[-,line width=2pt] (c) -- (f);
\draw[-,line width=2pt] (c) -- (g);
\draw[-,line width=1pt] (f) -- (h);
\draw[-,line width=1pt] (f) -- (i);
\draw[-,line width=2pt] (e) -- (j);
\draw[-,line width=2pt] (e) -- (k);
\end{tikzpicture}
\end{center}
\caption{A $2$-parking tree and 
the associated depth-first search ordering.
This example corresponds to the $2$-parking function $(3,1,9,1,4)$.}
\label{figure_p_tree}
\end{figure}
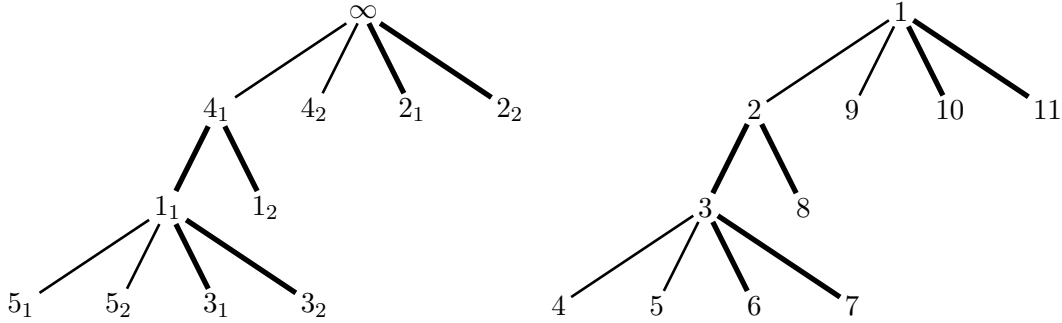

Note that for $d=1$ our definition yields a labeled rooted plane tree
with no restriction on the degrees of the vertices.  
A $2$-parking tree is shown on the left hand side of
Figure~\ref{figure_p_tree}. We consider $d$-parking trees as labeled
trees with {\em distinct vertices}. To identify the vertices in the
tree, in this paper we rely on the {\em depth-first search ordering
using the interval $[1,n]$}.
For a tree $T$ and a vertex $u$ in the tree~$T$, let $T(u)$ be the subtree
that has the node~$u$ as a root, that is, $T(u)$ consists of all
vertices that are descendants of $u$.

\begin{definition}
Given a plane tree $T$ on $n$ vertices and a positive integer $a$, the
{\em depth-first search ordering of $T$ using the interval $[a,a+n-1]$}
is defined by the following procedure:
\begin{itemize}
\item[(1)] We label the root $r$ of the tree $T$ with the number $a$,
that is, we set $\omega(r) = a$.
\item[(2)] Suppose $r_{1},r_{2},\ldots,r_{s}$ are the children of the root~$r$ in
  the left-to-right order and let $T_{i}$ be the subtree~$T(r_{i})$ for
  $i=1,2,\ldots,s$. We label each subtree $T_{i}$ recursively with
  depth-first search ordering using the
  interval $[a+|T_{1}|+\cdots+|T_{i-1}|+1,a+|T_{1}|+\cdots+|T_{i}|]$.
\end{itemize}
\end{definition}

The right hand side of Figure~\ref{figure_p_tree}
shows the {\em depth-first search ordering} of the vertices.
\begin{remark}
\label{remark_ourtrees}  
{\rm Similar $d$-parking trees (the mirror images of the
present ones and relying on the breadth-first search ordering) were defined
in~\cite[Section~6]{Hetyei}. Both possibilities
occur in the $d=1$ case as examples in~\cite[Section~13.2.3]{Yan_chapter}.
Another example of similar parking trees being
used may be found in~\cite{Ehrenborg_Hetyei_Readdy}, where the labels of
the tree on the left hand side appear as edge labels. The same idea of
shifting labels from edges to vertices and vice versa also appears
in~\cite{Irving_Rattan}. A key feature of all the 
parking trees discussed in~\cite[Section 13.2.3]{Yan_chapter} is that
$a_i=a_j$ holds in the parking function exactly when $i$ and $j$ label
siblings in the parking tree. This property is {\em not present} in the
parking trees of Kreweras~\cite{Kreweras_suites} and 
Knuth~\cite[Section~13.2.2]{Yan_chapter}, nor in the $d$-parking trees
of Irving and Rattan~\cite{Irving_Rattan}.}   
\end{remark}  

Recall there is a bijection between $d$-parking functions of length $k$ and
maximal chains in the noncrossing partition poset $\NC^{d}_{dk+1}$;
see Theorem~\ref{theorem_list_of_labels_bijection}.
We now include $d$-parking trees in this bijective family.
We begin with a key observation, that is applicable in all situations
when we identify the vertices of a plane tree using an {\em increasing ordering}.  
\begin{definition}
Given a plane tree $T$ on $n$ vertices an {\em increasing ordering} of the
vertices of $T$ is a bijection $\omega$ from the vertex set of $T$ to
the set $\{1,2,\ldots,n\}$ such that for every cover relation $u \prec v$ in the tree,
that is, $u$ is the parent of~$v$, the inequality $\omega(u) < \omega(v)$ holds.
\end{definition}  
Clearly, the depth-first search ordering and any breadth-first search
ordering are increasing orderings. It is worth noting that all {\em vertex
  orders} generated using {\em a choice function} in~\cite[Section
  13.2.3]{Yan_chapter} in recursive fashion have this property: the
interested reader may show this by induction on the number of vertices.    
The next lemma is a generalization of the observation made
in~\cite[Section~6.4, Exercise~31, solution (a)]{Knuth} stating the
following: given a plane tree on $n+1$ vertices, if
$(b_{1},\ldots,b_{n},0)$ is the list of numbers of children of its
vertices listed in the depth first search order, then any sequence
$(a_{1},\ldots,a_{n})$ containing $b_{i}$ copies of $i$ (for $1\leq i\leq
n$) is a parking function.  
\begin{lemma}
\label{lemma_parking}
Let $T$ be a $d$-parking tree on $n=dk+1$ vertices and let $\omega$ be
an increasing ordering of its vertices. For all $1 \leq i \leq k$
we define $a_{i}$ as the label $\omega(p)$ of the common parent of the
vertices $\omega(i_{1}),\omega(i_{2}), \ldots, \omega(i_{d})$.
Then the resulting vector $\vec{a}_{\omega}(T) = \vec{a} = (a_{1},a_{2}, \ldots, a_{k})$ is a
$d$-parking function. 
\end{lemma}  
\begin{proof}
Note that the number $j$ appears in the list $\vec{a}$ if and only if
the vertex $\omega^{-1}(j)$ is not a leaf. If $\omega^{-1}(j)$ has
$d \cdot c(j)$ children then $j$ appears in the list $\vec{a}$
exactly $c(j)$ times. Consider the ordered list
$a_{(1)} \leq a_{(2)} \leq \cdots \leq a_{(k)}$ obtained from $\vec{a}$. 
Let us relabel the vertices of the $d$-parking tree in such a way that
the order $(1_{1},1_{2},\ldots, 1_{d}, \ldots, k_{1}, k_{2},\ldots, k_{d})$
corresponds to listing all children of the vertex
$\omega^{-1}(a_{(1)})$, then of $\omega^{-1}(a_{(2)})$, and so on,
finally of $\omega^{-1}(a_{(k)})$  
in the left-to-right order. In the case when
$$a_{(i)}=a_{(i+1)}=\cdots =a_{(i+k)}=j$$
is a longest possible subsequence of repeated entries, we list
all children of $\omega^{-1}(j)$ at once in such a way that the
rules stated in Definition~\ref{definition_parking_tree} are observed.  

Let us call the resulting $d$-parking tree
the {\em straightened $d$-parking tree with respect to $\omega$}. The
straightened $2$-parking tree with respect to the depth-first search
order, obtained from the $2$-parking tree shown in
Figure~\ref{figure_p_tree} is represented in
Figure~\ref{figure_p_tree_straightened}. 
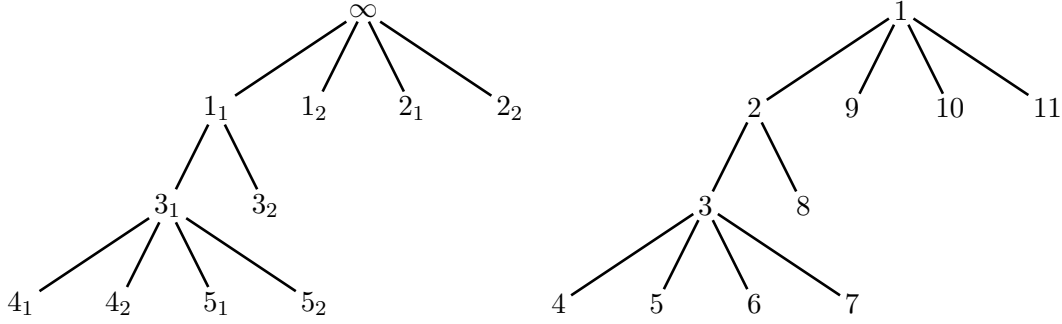
\begin{figure}[ht]
\newcommand{\shift}{5.5}
\begin{center}
\begin{tikzpicture}[scale = 1.3, inner sep=0.6mm]
\node (a) at (0,5) {$\infty$};
\node (b) at (-1.5,4) {$2_{1}$};
\node (c) at (-0.5,4) {$2_{2}$};
\node (d) at (0.5,4) {$1_{1}$};
\node (e) at (1.5,4) {$1_{2}$};
\node (f) at (-1,3) {$3_{1}$};
\node (g) at (0,3) {$3_{2}$};
\node (h) at (-1.5,2) {$4_{1}$};
\node (i) at (-0.5,2) {$4_{2}$};
\node (j) at (1,3) {$5_{1}$};
\node (k) at (2,3) {$5_{2}$};

\draw[-,line width=1pt] (a) -- (b);
\draw[-,line width=1pt] (a) -- (c);
\draw[-,line width=1pt] (a) -- (d);
\draw[-,line width=1pt] (a) -- (e);
\draw[-,line width=1pt] (c) -- (f);
\draw[-,line width=1pt] (c) -- (g);
\draw[-,line width=1pt] (f) -- (h);
\draw[-,line width=1pt] (f) -- (i);
\draw[-,line width=1pt] (e) -- (j);
\draw[-,line width=1pt] (e) -- (k);
\end{tikzpicture}
\hspace*{15mm}
\begin{tikzpicture}[scale = 1.3, inner sep=0.6mm]
\node (a) at (0,5) {$1$};
\node (b) at (-1.5,4) {$2$};
\node (c) at (-0.5,4) {$3$};
\node (d) at (0.5,4) {$8$};
\node (e) at (1.5,4) {$9$};
\node (f) at (-1,3) {$4$};
\node (g) at (0,3) {$7$};
\node (h) at (-1.5,2) {$5$};
\node (i) at (-0.5,2) {$6$};
\node (j) at (1,3) {$10$};
\node (k) at (2,3) {$11$};

\draw[-,line width=1pt] (a) -- (b);
\draw[-,line width=1pt] (a) -- (c);
\draw[-,line width=1pt] (a) -- (d);
\draw[-,line width=1pt] (a) -- (e);
\draw[-,line width=1pt] (c) -- (f);
\draw[-,line width=1pt] (c) -- (g);
\draw[-,line width=1pt] (f) -- (h);
\draw[-,line width=1pt] (f) -- (i);
\draw[-,line width=1pt] (e) -- (j);
\draw[-,line width=1pt] (e) -- (k);
\end{tikzpicture}
\end{center}
\caption{The straightened of the $2$-parking tree in Figure~\ref{figure_p_tree}
with respect to the depth-first search.
This example corresponds to the $2$-parking function $(1,1,3,4,9)$.}
\label{figure_p_tree_straightened}
\end{figure}
For $1 \leq i \leq k$ consider the list of children associated
to the parents in the initial segment $a_{(1)} \leq
a_{(2)} \leq \cdots \leq a_{(i-1)}$. In the straightened $d$-parking tree
this is exactly the list of the first $d\cdot(i-1)$ entries of the list $(1_{1},1_{2},\ldots, 1_{d}, \ldots, k_{1}, k_{2},\ldots,
k_{d})$. Observe that $a_{(i)}$ is exactly the value of
$\omega(p)$ for the common parent of $i_{1}, i_{2},\ldots,i_{d}$ in the
straightened parking tree.  Consider any nonroot-vertex $v$ that is not
yet listed. The parent $p(v)$ of $v$ satisfies $\omega(p(v))\geq
a_{(i)}$. Since the ordering $\omega$ is increasing on subtrees we must
also have $\omega(v)> a_{(i)}$. We obtained that the label of any
non-yet listed nonroot vertex is greater than the label $a_{(i)}$. There are
$n-(d(i-1)+1)$ such vertices which forces $a_{(i)} \leq d(i-1)+1$. 
\end{proof}  
As it was noted in~\cite[Section~1.6]{Irving_Rattan}, every parking
function gives rise to a $d$-parking function, obtained by replacing
each entry with $d$ copies of itself. In our model this expansion
operation gives rise to the following operation on our $d$-parking trees.
\begin{definition}
\label{definition_expansion}  
Let $T$ be a $d$-parking tree on $n=dk+1$ vertices. We call the {\em
  expansion of $T$} the $1$-parking tree obtained by replacing each
label $i_{j}$ with $(i-1)d+j$ for $1 \leq i \leq k$ and $1 \leq j \leq d$. 
\end{definition}  
The following statement is a direct consequence of
Lemma~\ref{lemma_parking} and Definition~\ref{definition_expansion}.
\begin{corollary}
Let $T$ be a $d$-parking tree on $n=dk+1$ vertices, let $\omega$ be
an increasing ordering of its vertices and let $\vec{a}_{\omega}(T)=\vec{a}=(a_{1},\ldots,a_{k})$ be
the $d$-parking function associated to $T$ ordered by~$\omega$ in Lemma~\ref{lemma_parking}.
Let $U$ be the expansion of the tree $T$.
Then the parking function $\vec{b}=\vec{a}_{\omega}(U)$
is the vector obtained  by replacing each entry $a_{i}$ in $\vec{a}$ by a
list of $d$ consecutive copies of $a_{i}$.
\end{corollary}  
\begin{example}
{\rm The $2$-parking tree shown in Figure~\ref{figure_p_tree}
corresponds to the $2$-parking function $\vec{a}=(3,1,9,1,4)$ with respect to he
depth-first search ordering. Its expansion corresponds to the parking
function $\vec{b}=(3,3,1,1,9,9,1,1,4,4)$. The fourth coordinate of
$\vec{a}$ is $1$ because the common parent $\infty$ of $4_{1}$ and $4_{2}$
satisfies $\omega(\infty)=1$. In the expansion $4_{1}$ becomes $3\cdot
2+1=7$ and $4_{2}$ becomes $3\cdot 2+2=8$. The seventh and eighth
coordinates of $\vec{b}$ are both $1$.}
\end{example}

In the remainder of the section we assume that
the increasing ordering $\omega$ is the depth-first search ordering.
In this case we will write 
$\vec{a}(T)$ for $\vec{a}_{\omega}(T)$
and call $\vec{a}(T)$
{\em the $d$-parking function associated to the $d$-parking tree $T$}.   

\begin{theorem}
The map $T \longmapsto \vec{a}(T)$ from $d$-parking trees on $n = dk+1$ vertices
to the $d$-parking functions of length $k$ is a bijection.
\label{theorem_parking_function_parking_tree}
\end{theorem}
\begin{proof}
It suffices to show that given a $d$-parking function $\vec{a} =
(a_{1},a_{2}, \ldots, a_{k})$ there is exactly one $d$-parking tree $T$
on $n=dk+1$ whose associated $d$-parking function is $\vec{a}$.

Observe first that we may restrict our attention to straightened
$d$-parking trees (defined in the proof of Lemma~\ref{lemma_parking})
and $d$-parking functions $\vec{a} = (a_{1},a_{2}, \ldots, a_{k})$
satisfying $a_{1} \leq a_{2} \leq \cdots \leq a_{k}$. Indeed, given any
permutation $\pi$ on the set $\{1,2,\ldots,k\}$, let $T_{\pi}$ be the
$d$-parking tree obtained from the $d$-parking tree $T$ by  replacing the labels
$i_{1},i_{2},\ldots, i_{d}$ with $\pi(i)_{1},\pi(i)_{2},\ldots,\pi(i)_{d}$ for each
$i\in \{1,2,\ldots,k\}$. If $\vec{a} = (a_{1},a_{2}, \ldots, a_{k})$ is
the $d$-parking function associated to $T$ then
$\vec{a}_{\pi} = (a_{\pi(1)},a_{\pi(2)}, \ldots, a_{\pi(k)})$ is the
$d$-parking function associated to $T_{\pi}$. In other words, once we
prove the statement for straightened $d$-parking trees and $d$-parking
functions $\vec{a} = (a_{1},a_{2}, \ldots, a_{k})$ satisfying $a_{1}\leq
a_{2}\leq \cdots \leq a_{k}$, the general statement may be obtained by
permuting the labels.

Consider now a $d$-parking function $\vec{a} = (a_{1},a_{2}, \ldots, a_{k})$
satisfying $a_{1} \leq a_{2} \leq \cdots \leq a_{k}$.
As in the proof of Lemma~\ref{lemma_parking} let $c(j)$ be the number of
times the number $j$ occurs in $\vec{a}$. Let
$\{j_{1}<j_{2}<\cdots<j_{m}\}$ be the set of integers satisfying $c(j) \geq 1$.
Since $a_{1}=1$, we must have $c(1)\geq 1$, $j_{1}=1$, and the
straightened $d$-parking tree $T$ giving rise to 
$\vec{a}$ must have exactly $d \cdot c(1)$ children: $1_{1},1_{2},\ldots, 1_{d},
\ldots, c(1)_{1}, c(1)_{2}, \ldots, c(1)_{d}$. The numbers $2,3,\ldots,j_{2}-1$
do not appear in $\vec{a}$, they must be numbers of leaves in the
depth-first search ordering, and these $j_{2}-2$ leaves must be the
leftmost $j_{2}-2$ children of the root. The vertex $j_{2}$ in the
depth-first search order is then child number $j_{2}-1$ of the root. This
assignment is possible if and only if
\begin{equation}
\label{inequality_parking_2}  
  j_{2} \leq d \cdot c(1) + 1
\end{equation}  
holds, as this is the number of children of the root. Observe also that we have
$a_{1}=a_{2}=\cdots=a_{c(1)}=1$ and $j_{2}$ is the common value of
$a_{c(1)+1}=a_{c(1)+2}=\cdots=a_{c(1)+c(j_{2})}$. The
inequality~\eqref{inequality_parking_2} is equivalent to
$a_{c(1)+1} \leq d \cdot c(1) + 1$ which is the $d$-parking condition for
$a_{c(1)+1}$. (It is worth noting that the parking conditions for
$a_{2}=\cdots=a_{c(1)}=1$  are automatically satisfied, and they follow
from~\eqref{inequality_parking_2} for
$a_{c(1)+2}=\cdots=a_{c(1)+c(j_{2})}$.)

We proceed now by induction of $m$, the number of distinct values listed
in $\vec{a}$. Assume that we have already found a unique straightened
parking tree $T'$ on $d \cdot (c(j_{1})+c(j_{2})+\cdots+c(j_{m-1}))+1$ vertices whose
associated $d$-parking function is
$(a_{1},a_{2},\ldots,a_{c(j_{1})+c(j_{2})+\cdots+c(j_{m-1})})$. By our
construction, $v$ is a nonleaf vertex of $T'$ if and only if $\omega(v)$
belongs to the set $\{j_{1},j_{2},\ldots,j_{m-1}\}$.  Now we want to
add $d \cdot c(j_{m})$ children to the vertex labeled $j_{m}$ in $T'$, to obtain
the tree $T$. In the depth-first search ordering this will increase the
value of $\omega$ by $d \cdot c(j_{m})$ for the vertices $v$ of $T$ that satisfy
$\omega(v)>j_{m}$, but by $j_{m-1}<j_{m}$ all these vertices are leaves, and
the $d$-parking function associated to $T$ will agree with the one
associated to $T'$ in the first $c(j_{1})+c(j_{2})+\cdots+c(j_{m-1})$
  coordinates. The insertion of children is possible if and only if
\begin{equation}
\label{inequality_parking_m}  
j_{m-1} < j_{m} \leq d \cdot (c(j_{1})+c(j_{2})+\cdots+c(j_{m-1})) + 1
\end{equation}  
holds. The upper bound is equivalent to
\begin{align*}
a_{c(j_{1})+c(j_{2})+\cdots+c(j_{m-1})+1}
& \leq d \cdot (c(j_{1})+c(j_{2})+\cdots+c(j_{m-1})) + 1,
\end{align*}
the $d$-parking condition for
$a_{c(j_{1})+c(j_{2})+\cdots+c(j_{m-1})+1}$. If this condition is satisfied,
the weaker conditions for $a_{c(j_{1})+c(j_{2})+\cdots+c(j_{m-1})+1},
\ldots, a_{c(j_{1})+c(j_{2})+\cdots+c(j_{m})}$ are also satisfied. The
inequality
\begin{align*}
j_{m-1} < d \cdot (c(j_{1})+c(j_{2})+\cdots+c(j_{m-1})) + 1
\end{align*}
needed to have at least one choice for $j_{m}$ is also a direct
consequence of the $d$-parking condition for
$a_{c(j_{1})+c(j_{2})+\cdots+c(j_{m-2})+1}$.   
\end{proof}

Theorem~\ref{theorem_list_of_labels_bijection}
establishes a bijection between the set of all maximal chains in
$\NC^{d}_{dk+1}$ and the set of all $d$-parking functions of length~$k$.
Composing this bijection with the one stated in
Theorem~\ref{theorem_parking_function_parking_tree} above we obtain a
bijection between the set of all maximal chains in
$\NC^{d}_{dk+1}$ and the set of all $d$-parking trees on $n=dk+1$
vertices. The correspondence is visually straightforward, as exhibited
in the following example.

\begin{example}
\label{example_chain}  
{\rm
Consider the $2$-parking function $(a_{1}, \ldots, a_{5}) = (3,1,9,1,4)$,
corresponding to the $2$-parking tree in Figure~\ref{figure_p_tree}.
The common parent of the nodes labeled $5_{1}$ and $5_{2}$ is
labeled $1_{1}$, and the depth-first search label of $1_{1}$ is
$\omega(1_{1}) = 4 = a_{5}$. This $2$-parking function 
corresponds to the following maximal chain in the noncrossing partition
poset $\NC^{2}_{11}$: 
\begin{align*}
\hz
& \coveredby
1|2|3,4,7|5|6|7|8|9|10|11
\coveredby
1,8,9|2|3,4,7|5|6|10|11\\
& \coveredby
1,8,9,10,11|2|3,4,7|5|6|9
\coveredby
1,2,3,4,7,8,9,10,11|5|6
\coveredby
\ho .
\end{align*}
In the $4$th step we join together the blocks containing the elements
$\omega(\infty) = 1$,
$\omega(4_{1}) = 2$ and
$\omega(4_{2}) = 3$.
Here $\infty$ is the common parent of $4_{1}$ and $4_{2}$. The
edges of the nonsingleton blocks $\{1,8,9,10,11\}$ and $\{3,4,7\}$ are
shown in bold in the figure. 
It is worth noting that in the edge labeling introduced
in~\cite{Muhle_Nadeau_Williams} the label of this cover relation is $1$
because that is the least element of the cycle $(1,2,7)$ and
$$
(3,4,7)(1,8,9,10,11)(1,2,7)=(1,2,3,4,8,9,10,11). 
$$
There seems to be no obvious rule to identify the cycle $(1,2,7)$ in the
picture, but the blocks that are merged are clearly recognizable. 
}
\end{example}
In general, in step $i$ we join the blocks
containing the elements $\omega(i_1),\omega(i_2),\ldots,\omega(i_d)$
with the block containing the element $\omega(p(i))$. If we order the
blocks in increasing order of their minimum elements, then $B_1$ is the
block containing $\omega(p(i))$, and for $1\leq s\leq d$ we have $\min
(B_{s+1})=\omega(i_{s})$, since we have
$$\omega(p(i))<\omega(i_1)<\omega(i_2)<\cdots< \omega(i_d).$$
It is easy to see that $\omega(p(i))$ is the largest element of $B_1$ that
is less than $\min(B_2)$. Indeed, any element of $B_1$ that succeeds
$\omega(p(i))$ in the depth-first search order is at least as large as
$\omega(i_1)=\min(B_2)$, the leftmost child of $\omega(p(i))$. We only
need to show that the blocks arising in this process are indeed the blocks of
the noncrossing partitions that occur in the corresponding maximal
chain. This is the result stated in
Theorem~\ref{theorem_parking_noncrossing} below. To prove it, we will
need the next lemma and the following procedure. 
The lemma highlights details regarding the
bijection in Theorem~\ref{theorem_list_of_labels_bijection}
and the procedure shows how to compute the inverse of this bijection.
 
\begin{lemma}
Let ${\bf m} = \{\hz = \pi_{0} \coveredby \pi_{1} \coveredby \cdots \coveredby \pi_{k} = \ho\}$
be a maximal chain in $\NC^{d}_{dk+1}$
and
let $\lambda({\bf m}) = (a_{1},a_{2}, \ldots, a_{k})$.
Let $r$ be the largest label occurring among these labels,
that is, $r = \max(a_{1}, a_{2}, \ldots, a_{k})$.
Furthermore, let $s$ be the last position this label occurs, that is,
$s = \max(\{i : a_{i} = r\})$. 
Then the partition $\pi_{s-1}$ contains the singleton blocks
$\{r+1\}, \{r+2\}, \ldots, \{r+d\}$.
\label{lemma_d_singleton_blocks}
\end{lemma}
\begin{proof}
We claim that for $0 \leq e \leq d$
the partition $\pi_{s-1}$ contains the singleton blocks
$\{r+1\}$, $\{r+2\}$, $\ldots$, $\{r+e\}$.
We prove this by induction on $e$.
The basis case is $e=0$ which is directly true.

Assume now that the statement is true for $0 \leq e \leq d-1$ and
we prove it for $e+1$.
Let $B$ be the block of $\pi_{s-1}$ that contains the element $r+e+1$.
We claim that the block $B$ does not contain the element $r$.
If $e=0$ note that $\{r,r+1\} \subseteq B$ contradicts that
$\lambda(\pi_{s-1},\pi_{s}) = r$.
If $e \geq 1$ 
then we have a block in the dual partition whose cardinality lies strictly
between $1$ and $d+1$, which leads to a contradiction.
Hence we may assume that $r$ and $r+e+1$ lie in different blocks.
If $r+e+1$ is the smallest element in the block~$B$, then this
contradicts that the largest label of the maximal chain is~$r$.
Hence assume there is an element $t$ such that $1 \leq t < r$ in the block~$B$.
Let $B_{1}$ be the block of~$\pi_{s-1}$ that contains the element $r$.
Since $\lambda(\pi_{s-1},\pi_{s}) = r$ we know that $B_{1}$ is
joined with $d$ other blocks $B_{2}$ through $B_{d+1}$
to obtain the partition $\pi_{s}$.
These other blocks must have their minimal elements greater than $r$.
By the noncrossing property, the only such blocks are
the singletons $\{r+1\}$ through $\{r+e\}$. But there are fewer than $d$
of these blocks. Hence there is no such element $t$ and we conclude
that the element $r+e+1$ forms a singleton block, proving the induction step.
\end{proof}

\begin{procedure}
Given a $d$-parking function $\vec{a} = (a_{1}, a_{2}, \ldots, a_{k})$
we can reconstruct the maximal chain ${\bf m}$ in $\NC^{d}_{dk+1}$ such that
$\lambda({\bf m}) = \vec{a}$ by the following recursive procedure:
If $k=0$ let ${\bf m}$ be the unique maximal chain in $\NC^{d}_{1}$.
For $k \geq 1$
let $r =\max(a_{1}, a_{2}, \ldots, a_{k})$ and $s = \max(\{i : a_{i} = r\})$
as in Lemma~\ref{lemma_d_singleton_blocks}.
Now
\begin{align*}
\vec{b} & = (b_{1}, b_{2}, \ldots, b_{k-1}) = (a_{1}, a_{2}, \ldots, a_{s-1}, a_{s+1}, \ldots, a_{k})
\end{align*}
is a $d$-parking function.
By recursion we have a maximal chain
${\bf n} = \{\hz = \sigma_{0} \coveredby \sigma_{1} \coveredby \cdots \coveredby \sigma_{k-1} = \ho\}$
in $\NC^{d}_{d(k-1)+1}$
such that
$\lambda({\bf n}) = \vec{b}$.
Let $f : [d(k-1)+1] \longrightarrow[dk+1]$
be the following relabeling function: we set
\begin{align}
f(j)
& =
\begin{cases}
j     & \text{ if } 1    \leq j \leq r, \\
j+d & \text{ if } r+1 \leq j \leq d(k-1)+1.
\end{cases}
\label{equation_relabeling_function}
\end{align}
Define the maximal chain
${\bf m} = \{\hz = \pi_{0} \coveredby \pi_{1} \coveredby \cdots \coveredby \pi_{k} = \ho\}$
by 
\begin{align*}
\pi_{i}
& =
\begin{cases}
\{f(B) \: : \: B \in \sigma_{i}\} \cup \{\{r+1\}, \ldots, \{r+d\}\}
& \text{ for } 0 \leq i \leq s-1, \\
\{f(B) \: : \: r \notin B \in \sigma_{i-1}\}
\cup
\{f(B) \cup \{r+1, \ldots, r+d\} \: : \: r \in B \in \sigma_{i-1}\}
& \text{ for } s \leq i \leq k.
\end{cases}
\end{align*}
\label{procedure_reconstructing_maximal_chain}
\end{procedure}
\begin{proof}
Note that the relabeling function $f$ opens up for an interval of $d$ new elements,
that is, $r+1$ through $r+d$.
Let $a$ and $b$ be two elements in the cyclic order on $[d(k-1)+1]$.
Then the number elements between $a$ and $b$ in the cyclic order is congruent to the
number of elements between $f(a)$ and~$f(b)$ in the cyclic order on
$[dk+1]$ modulo $d$.
Hence by Theorem~\ref{theorem_characterization}
we know that the partition~$\pi_{i}$ belongs to the poset $\NC^{d}_{dk+1}$ for all $i$.
We also obtain $\lambda({\bf m}) = \vec{a}$.
\end{proof}

\begin{theorem}
Let ${\bf m} = \{\hz = \pi_{0} \coveredby \pi_{1} \coveredby \cdots \coveredby \pi_{k} = \ho\}$
be a maximal chain in $\NC^{d}_{dk+1}$
and let $T$ be the $d$-parking tree such that $\lambda({\bf m}) = \vec{a} = \vec{a}(T)$. 
Then $\pi_{i}$ is the noncrossing partition whose blocks consist of the connected
components of the graph containing all edges that connect a vertex~$j_{s}$
to its parent for some $1 \leq j \leq i$ and some $1 \leq s \leq d$.
In particular, the cover relation from the partition $\pi_{i-1}$ to the
partition $\pi_{i}$ is obtained by joining the blocks containing the $d+1$ elements
$\omega(p_{i}), \omega(i_{1}), \omega(i_{2}), \ldots, \omega(i_{d})$, where
$p_{i}$ is the common parent of $i_{1}$ through $i_{d}$. 
\label{theorem_parking_noncrossing}
\end{theorem}
\begin{proof}
We proceed by induction on $k$.
When $k=0$ there is is nothing to prove completing the induction basis.
Assume now that the result is true for $k-1$ and we prove it for $k$.

Let $r$ and $s$ be defined
as in Lemma~\ref{lemma_d_singleton_blocks}.
That is, $r$ is the largest entry in the $k$-parking function $\vec{a} = (a_{1},a_{2}, \ldots, a_{k})$
and $s$ is the last entry such that $a_{s} = r$.
Lemma~\ref{lemma_d_singleton_blocks} implies that
the partition $\pi_{s-1}$ contains the singleton blocks
$\{r+1\}, \{r+2\}, \ldots, \{r+d\}$.
As in Procedure~\ref{procedure_reconstructing_maximal_chain}
let $\vec{b} = (b_{1}, b_{2}, \ldots, b_{k-1}) = (a_{1}, a_{2}, \ldots, a_{s-1}, a_{s+1}, \ldots, a_{k})$
be the $d$-parking function where we remove the $a_{s}$ entry.
Let $U$ be the $d$-parking tree corresponding to the $d$-parking function~$\vec{b}$,
that is, we define $U$ by $\vec{a}(U) = \vec{b}$.

Construct a new tree $T^{*}$ by the following two steps:
\begin{itemize}
\item[(1)]
Relabel the nodes $i_{j}$ in the tree $U$ where $i \geq s$ to be $(i+1)_{j}$.
Note that after the relabeling, there are no nodes with the labels $s_{1}$ through $s_{d}$.
\item[(2)]
Let $p$ be the node in the tree $U$ such that $\omega_{U}(p) = r$.
Attach the new leaves $s_{1}$ through $s_{d}$ to the node $p$
such that they are the $d$ right-most children of the node $p$.
Note that these nodes in the depth-first search labeling
receive the labels $r+1$ through $r+d$, that is,
$\omega_{T^{*}}(s_{j}) = r+j$ for $1 \leq j \leq d$.
\end{itemize}
The node $x$ in the tree $U$ satisfies
$f(\omega_{U}(x)) = \omega_{T^{*}}(x)$ where $f$ is the relabeling function
defined in~\eqref{equation_relabeling_function}
Hence $\vec{a}(T^{*}) = \vec{a}$ holds.
Since the correspondence between $d$-parking trees and $d$-parking functions
is a bijection, we conclude that $T^{*}$ is the tree $T$.

Consider the cover relation $\pi_{i-1} \coveredby \pi_{i}$
in the chain ${\bf m}$.
By the induction hypothesis
when $i \neq s$
the partition $\pi_{i}$ is obtained from $\pi_{i-1}$ by joining the blocks containing the elements
$f(\omega_{U}(p_{i}))$, $f(\omega_{U}(i_{1}))$, $f(\omega_{U}(i_{2}))$, $\ldots$, $f(\omega_{U}(i_{d}))$.
But these are exactly the elements
$\omega_{T}(p_{i}), \omega_{T}(i_{1}), \omega_{T}(i_{2}), \ldots, \omega_{T}(i_{d})$.
Finally, when $i$ is $s$ by the above construction we are joining
$\omega_{T}(p), \omega_{T}(s_{1}), \omega_{T}(s_{2}), \ldots, \omega_{T}(s_{d})$
which are the elements $r$ through $r+d$, completing the induction.
\end{proof}

It is an immediate consequence of the definitions that for a $d$-parking
function $\vec{a}=(a_{1},a_{2},\ldots,a_{k})$ 
and a permutation $\tau$ of the set $\{1,2,\ldots,k\}$ the vector
$\tau(\vec{a})=(a_{\tau(1)},a_{\tau(2)},\ldots,a_{\tau(k)})$ is also a
$d$-parking function. Since $d$-parking functions bijectively label the
maximal chains of $\NC^{d}_{dk+1}$, the above action of the symmetric
group $\mathfrak{S}_{k}$ on the set of $d$-parking functions of
length $k$ induces an action of the same group on the maximal chains of
$\NC^{d}_{dk+1}$. This action is described
in~\cite[Lemma~5.2]{Muhle_Nadeau_Williams} and used in the proof
of~\cite[Theorem~5.3]{Muhle_Nadeau_Williams}. As a direct consequence of
Theorem~\ref{theorem_parking_noncrossing} we obtain the following
statement.

\begin{corollary}
\label{corollary_permutations}  
Two maximal chains of $\NC^{d}_{dk+1}$ are in the same orbit of the
above described action of the symmetric group if and only if the
removal of the labels on the $d$-parking trees associated to them yields
the same plane tree. 
\end{corollary}   

Theorem~\ref{theorem_parking_noncrossing} and
Corollary~\ref{corollary_permutations} have the following consequence.

\begin{corollary}
\label{corollary_parkingblocks}
Let $n=dk+1$. If a partition $\pi$ of $[n]$ belongs to the noncrossing
partition poset $\NC^{d}_{n}$ then there is a $d$-parking tree
$T$ and an $s\in\{0,1,\ldots,k\}$ such that the blocks of $\pi$
are the connected components of the graph obtained from the tree $T$ by
deleting all the edges connecting each vertex $i_{j}$ to its parent
for all $i>s$ and $j\in\{1,2,\ldots,d\}$. Conversely, if there is a
$d$-parking tree $T$ and a subset $S$ of $\{1,2,\ldots,k\}$, such that
the connected components of the graph obtained from the tree $T$ by
deleting all the edges connecting each vertex $i_{j}$ to its parent
for all $i\in S$ and $j\in\{1,2,\ldots,d\}$ are the blocks of~$\pi$, then
$\pi$ belongs to $\NC^{d}_{n}$.
\end{corollary}  

For $d=1$ Corollary~\ref{corollary_parkingblocks} yields the following
statement. 

\begin{corollary}
\label{corollary_ncparkingblocks}
A partition $\pi$ of $[n]$ is a noncrossing partition
if and only if its blocks are the
connected components of a graph obtained from a plane tree ordered by
the depth-first search order after deleting an arbitrary subset of its
edges.
\end{corollary}

\section{Concluding remarks}
\label{section_Concluding_remarks}

As it was first pointed out
in~\cite{Muhle_Nadeau_Williams}, many of the enumerative results are
related to the Raney numbers.
Some of these results already have a combinatorial proof.
For instance, is there a bijective proof of Theorem~\ref{theorem_1_d+1}?
Could the tree techniques of Bernardi and Morales~\cite{Bernardi_Morales}
be useful?
A slightly different question, 
is there a combinatorial proof for 
Corollary~\ref{corollary_falling_chains}?

The quasisymmetric function of a graded poset
encodes all of the flag $f$-vector information of the poset;
see~\cite{Ehrenborg}.
Stanley observed that if every interval of a poset $P$
is self-dual then the quasisymmetric function of $P$
is a symmetric function; see~\cite[Theorem~1.4]{Stanley_Flag}.
Hence the noncrossing partition poset~$\NC^{d}_{n}$
has a symmetric quasisymmetric function.
In the paper~\cite{Stanley}
Stanley explores the quasisymmetric function
of the noncrossing partition lattice $\NC_{n}$.
Are there similar results for the $d$-indivisible noncrossing partition
poset $\NC^{d}_{n}$? 

The set of $d$-parking functions are known to be in bijection with
the regions of the extended Shi-arrangement via the Pak-Stanley
labeling~\cite[Theorem~2.1]{Stanley_Flag}. The same
$d$-parking functions also label the maximal chains of the poset
$\NC^{d}_{n}$. Is there a geometric way to directly connect the poset
$\NC^{d}_{n}$ with the extended Shi arrangement?

Noncrossing partitions have connections with free probability.
Does the subposet $\NC^{d}_{n}$ of the lattice $\NC_{n}$ have a similar
connection? 

There are noncrossing partitions for other Coxeter systems;
see~\cite{Armstrong,Reiner}, a type~$B$ analogue of the $d$-indivisible
noncrossing partitions has been proposed by M\"uhle, Nadeau and
Williams~\cite{Muhle_Nadeau_Williams}.  
The authors are currently developing the analogous results for
noncrossing partitions of these types.

As noted in Section~\ref{section_tree_representation}, a more complex
bijection between labeled trees and $k$-parking functions was constructed by
Irving and Rattan~\cite{Irving_Rattan}. Their bijection provided a
bijective proof of Yan's equidistribution
result~\cite{Yan}. It would be interesting to
find out what analogue of this result (if any) could be shown using the
bijection described in Section~\ref{section_tree_representation}. The
special case when $d=1$ is already discussed in~\cite[Section~13.2.3]{Yan_chapter},
but even that discussion may be worth a second look. 

We found the $EL$-labeling of $\NC^{d}_{n}$ by generalizing Stanley's
approach~\cite{Stanley} establishing a direct link between parking
functions and maximal chains of the noncrossing partition lattice. The
approach taken by M\"uhle, Nadeau and
Williams~\cite{Muhle_Nadeau_Williams} is closer to that of
Brady~\cite{Brady}  who established a direct bijection between minimal
factorizations of the full cycle $(1,2,\ldots,n)$ in the symmetric group
$\mathfrak{S}_{n}$ into transpositions and maximal chains in the noncrossing
partition lattice on the set $[n]$. It would be worth exploring further
whether the two approaches could be brought closer to each other by
studying further the parking trees we introduced in
Section~\ref{section_tree_representation}. 

\section*{Acknowledgements}

We are grateful to two anonymous referees for their helpful comments and
suggestions on an earlier draft of this paper. We thank Natasha Blitvi\'c
for bringing Beissinger's work to our attention. This work was partially
supported by a grant from the Simons Foundation
(\#854548 to Richard Ehrenborg and \#514648 to G\'abor Hetyei).

\newcommand{\journal}[6]{{\sc #1,} #2, {\it #3} {\bf #4} (#5), #6.}
\newcommand{\book}[4]{{\sc #1,} ``#2,'' #3, #4.}
\newcommand{\bookf}[5]{{\sc #1,} ``#2,'' #3, #4, #5.}
\newcommand{\arxiv}[3]{{\sc #1,} #2, {\tt #3}.}
\newcommand{\thesis}[4]{{\sc #1,} ``#2,'' Doctoral dissertation, #3,~#4.}
\newcommand{\preprint}[3]{{\sc #1,} #2, preprint {(#3)}.}
\newcommand{\preparation}[2]{{\sc #1,} #2, in preparation.}
\newcommand{\toappear}[3]{{\sc #1,} #2, to appear in {\it #3}.}

\end{document}